\documentclass{article}
\usepackage{amssymb}

\usepackage{graphicx}
\usepackage{amsmath}

\newtheorem{theorem}{Theorem}

\newtheorem{definition}[theorem]{Definition}

\input{tcilatex}

\begin{document}

\title{The heritage of S. Lie and F. Klein: \\
Geometry via transformation groups}
\author{Erc\"{u}ment Orta\c{c}gil}
\maketitle

\begin{abstract}
We outline a framework which generalizes Felix Klein's \textit{Erlangen
Programm }which he announced in 1872 after exchanging ideas with Sophus Lie%
\textit{.}
\end{abstract}

\section{$\protect\bigskip $Introduction}

Let $G$ be a Lie group which acts transitively on some space $M.$ In this
framework, Felix Klein defined geometry (later known as \textit{Erlangen
Programm}) as the study of properties which are invariant under the action
of $G.$ Realizing $M$ as the coset space $G/H$ where $H\subset G$ is the
stabilizer of some point $p\in M$, we may thus speak of the Klein pair (or
Klein geometry) $(G,G/H)$. In Euclidean geometry, the best known example, $G$
is the isometry group of $\mathbb{R}^{n}$ and is the semi-direct product of
the orthogonal group $O(n)$ and $\mathbb{R}^{n},$ $H=O(n)$ and $G/H=\mathbb{R%
}^{n}.$ Non-Euclidean geometries, which were known at Klein's time, are
other examples. Riemannian geometry can not be realized in this way unless
the metric has constant curvature in which case the isometry group acts
transitively and we have again a Klein geometry. We will refer the reader to
[14], pg. 133 so that he/she can feel the philosophical disturbance created
by this situation at that time. As an attempt to unify (among others)
Riemannian geometry and \textit{Erlangen Programm}, Elie Cartan introduced
in 1922 his generalized spaces (principal bundles) which are curved analogs
of the principal $H$-bundle $G\rightarrow G/H.$ We will mention here the two
outstanding books: [29] for Klein geometries, their generalizations called
Cartan geometries in [29] and also for the notions of Cartan and Ehresmann
connections, and [14] for more history of this subject (see in particular
pg. 34-42 for \textit{Erlangen Programm)}. Cartan's approach, which is later
devoloped mainly by topologists from the point of view of fiber bundles,
turned out to be extremely powerful. The spectacular achievements of the
theory of principle bundles and connections in geometry, topology and
physics are well known and it is therefore redundant to elobarate on them
here. However, it is also a fact that this theory leaves us with the
unpleasent question: What happened to \textit{Erlangen Programm}? The main
reason for this question is that the total space $P$ of the principal bundle 
$P\rightarrow M$ does not have any group-like structure and therefore does
not act on the base manifold $M.$ Thus $P$ emerges in this framework as a
new entity whose relation to the geometry of $M$ may not be immediate and we
must deal now with $P$ as a seperate problem. Consequently, it seems that
the most essential feature of Klein's realization of geometry is given up by
this approach. Some mathematicians already expressed their dissatisfaction
of this state of affairs in literature with varying tones, among which we
will mention [29], [35], [26] and other contraversial works of the author of
[26] (see the references in [26]).

The purpose of this work is to present another such unification which we by
no means claim to be the ultimate and correct one but believe that it is
faithful to Klein's original conception of geometry. This unification is
based on the ideas which S. Lie and F. Klein communicated to each other
before 1872 (see [14] for the extremely interesting history of this subject)
and the works of D.C. Spencer and his co-workers around 1970 on the formal
integrability of PDEs. The main idea is simple and seems to be the only
possible one: We concentrate on the action of $G$ on $M=G/H$ and generalize
this action rather than generalizing the action of $H$ on $G$ in the
principal $H$-bundle $G\rightarrow G/H.$ Now $G\subset Diff(M)$ and $G$ may
be far from being a Lie group. As the natural generalization of \textit{%
Erlangen Programm, }we may now deal directly with the group $G$ as in [2]
(see also [23]), but this approach again does not incorporate Riemannian
geometry unless the metric is homogeneous. We consider here the Lie
pseudogroup $\widetilde{G}$ determined by $G$ and filter the action of $%
\widetilde{G}$ on $M$ via its jets, thus realizing $\widetilde{G}$ as a
projective limit $Lim_{\leftarrow k}$ $\mathcal{S}_{k}(M)$ of Lie equations $%
\mathcal{S}_{k}(M).$ Lie equations (in finite form) are very special
groupoids and are very concrete objects which are extensively studied by
Spencer and his co-workers culminating in the difficult work [11]. We will
refer to [11], [25], [26] for Lie equations and [19], [20] for
differentiable groupoids and algebroids. On the infinitesimal level, we
obtain the approximation $Lim_{\leftarrow k}$ $\frak{s}_{k}(M)$ where $\frak{%
s}_{k}(M)$ is the infinitesimal Lie equation (or the algebroid) of $\mathcal{%
S}_{k}(M).$ The idea is now to start with the expression $Lim_{\leftarrow k}$
$\mathcal{S}_{k}(M)$ as our definition of homogeneous geometry $\mathcal{S}%
_{\infty }(M)$ (Section 3, Definition 2). Any transitive pseudogroup (in
particular a complex, symplectic or contact structure) determines a
homogeneous geometry and Klein geometries are special cases (Section 4).
Some $\mathcal{S}_{k}(M)$ may not prolong to a homogeneous geometry due to
the lack of formal integrability and Riemannian geometry (almost complex,
almost symplectic...structures) emerge as truncated geometries (Section 5).
We associate various spectral sequences to a homogeneous geometry $\mathcal{S%
}_{\infty }(M)$ (in particular to a truncated geometry $\mathcal{S}_{k}(M))$
(Sections 2, 3). For a complex structure, we believe that one these spectral
sequences is directly related to the Fr\"{o}licher spectral sequence which
converges to de Rham cohomology with $E_{1}$ term equal to Dolbeault
cohomology. This unification is also a natural specialization of the
standard principal bundle approach initiated by E. Cartan (Sections 6, 7, 8).

The idea of filtering an object via jets (for instance, the solution space
of some PDE called a diffiety in [32], [33]) is not new and is used by
A.M.Vinogradov in 1978 in the construction of his $\mathcal{C}$-spectral
sequence and in the variational bicomplex approach to Euler-Lagrange
equations (see [32], [33] and the references therein). In fact, this paper
can also be considered as the first step towards the realization of a
program stated in [33] for quasi-homogeneous geometries ([33], Section 6.4).
Further, there is definitely a relation between the higher order de Rham
complexes constructed here and those in [31]. We also believe that the main
idea of this paper, though it may not have been stated as explicitly as in
this paper, is contained in [26] and traces back to [11] and [25]. In
particular, we would like to emphasize that all the ingredients of this
unification are known and exist in the literature.

This paper consists of nine sections. Section 2 contains the technical core
in terms of which the geometric concepts substantiate. This section may be
somewhat demanding for the reader who is not much familiar with jets and the
formalism of Spencer operator. However, as we proceed, technical points
slowly evaporate and the main geometric concepts which we are all familiar
with, start to take the center stage.

\section{Universal homogeneous envelope}

Let $M$ be a differentiable manifold and $Diff(M)$ be the group of
diffeomorphisms of $M.$ Consider the map $Diff(M)\times M\rightarrow M$
defined by $(g,x)\rightarrow g(x)=y$ and let $j_{k}(g)_{y}^{x}$ denote the $%
k $-jet of $g$ with source at $x$ and target at $y.$ By choosing coordinates 
$(U,x^{i})$ and $(V,y^{i})$ around the points $x,y,$ we may think $%
j_{k}(g)_{y}^{x}$ as the coefficients of the Taylor expansion of $g(x)=y$ up
to order $k.$ We will call $j_{k}(g)_{y}^{x}$ the $k$-arrow induced by $g$
with source at $x$ and target at $y$ and imagine $j_{k}(g)_{y}^{x}$ as an
arrow starting at $x$ and ending at $y.$ Let $(f_{k})_{y}^{x}$ denote 
\textit{any }$k$-arrow, i.e., $(f_{k})_{y}^{x}$ is the $k$-jet induced by
some arbitrary local diffeomorphism which maps $x$ to $y$. With some
assumptions on orientability which involve only 1-jets (see [24] for
details), there exists some $g\in Diff(M)$ with $%
j_{k}(g)_{y}^{x}=(f_{k})_{y}^{x}.$ Therefore, with the assumption imposed by
[24], the pseudogroup $\widetilde{Diff(M)}$ of local diffeomorphisms on $M$
induces the same $k$-arrows as $Diff(M),$ but this fact will not be used in
this paper. We can compose successive $k$-arrows and invert all $k$-arrows.

Now let $(\mathcal{G}_{k})_{y}^{x}$ denote the set of all $k$-arrows with
source at $x$ and target at $y.$ We define $\mathcal{G}_{k}(M)\doteq \cup
_{x,y\in M}(\mathcal{G}_{k})_{y}^{x}$ and obtain the projections $\pi _{k,m}:%
\mathcal{G}_{k}(M)\rightarrow \mathcal{G}_{m}(M),$ $1\leq m\leq k-1,$ and $%
\pi _{k,m}$ is compatible with composition and inversion of arrows. We will
denote all projections arising from the projection of jets by the same
notation $\pi _{k,m}.$ Now $\mathcal{G}_{k}(M)$ is a transitive Lie equation
in finite form $(TLEFF)$ on $M$ which is a very special groupoid (see [11],
[25], [26] for Lie equations and [19], [20] for groupoids). We also have the
locally trivial map $\mathcal{G}_{k}(M)\rightarrow M\times M$ which maps $(%
\mathcal{G}_{k})_{y}^{x}$ to $(x,y).$ Note that $(\mathcal{G}_{k})_{x}^{x}$
is a Lie group and can be identified (not in a canonical way) with $k^{th}$%
-order jet group. Thus we obtain the sequence of \ homomorphisms

\begin{equation}
......\longrightarrow \mathcal{G}_{k+1}(M)\longrightarrow \mathcal{G}%
_{k}(M)\longrightarrow .....\longrightarrow \mathcal{G}_{1}(M)%
\longrightarrow M\times M\longrightarrow 1
\end{equation}

where the last arrow is used with no algebraic meaning but to express
surjectivity. (1) gives the vague formula $Diff(M)\times M=Lim_{k\rightarrow
\infty }\mathcal{G}_{k}(M)$ or more precisely $\widetilde{Diff(M)}%
=Lim_{k\rightarrow \infty }\mathcal{G}_{k}(M).$ The ambiguity in this last
formula is that a formal Taylor expansion may fail to recapture a local
diffeomorphism. However, this ambiguity is immaterial for our purpose for
the following reason: Let $(j_{\infty }g)_{x}^{x}$ denote the $\infty $-jet
of some local diffeomorphism $g$ where $g(x)=x.$ Now $(j_{\infty }g)_{x}^{x}$
determines $g$ modulo the $\infty $-jet of the identity diffeomorphism. This
is a consequence of the following elementary but remarkable fact: For 
\textit{any }sequence of real numbers $a_{0},a_{1},....$, there exists a
real valued differentiable function $f$ defined, say, near the origin $o\in 
\mathbb{R}$, satisfying $f^{(n)}(o)=a_{n}.$ In particular, the same
interpretation is valid for the $\infty $-jets of all objects to be defined
below.

Since $\mathcal{G}_{k}(M)$ is a differentiable groupoid (we will call the
object called a Lie groupoid in [19], [20] a differentiable groupoid,
reserving the term ``Lie'' for Lie equations), it has an algebroid $\frak{g}%
_{k}(M)$ which can be constructed using jets only. To do this, we start by
letting $J_{k}(T(M))_{p}$ denote the vector space of $k$-jets of vector
fields at $p\in M$ where $T(M)\rightarrow M$ is the tangent bundle of $M.$
An element of $J_{k}(T(M))_{p}$ is of the form $(p,\xi ^{i}(p),\xi
_{j_{1}}^{i}(p),\xi _{j_{2}j_{1}}^{i}(p),....,\xi
_{j_{k}j_{k-1}....j_{1}}^{i}(p))$ in some coordinates $(U,x^{i})$ around $p.$
If $X=(\xi ^{i}(x))$, $Y=(\eta ^{i}(x))$ are two vector fields on $U,$
differentiating the usual bracket formula $[X,Y](x)=\xi ^{a}(x)\partial
_{a}\eta ^{i}(x)-\eta ^{a}(x)\partial _{a}\xi ^{i}(x)$ successively $k$%
-times and evaluating at $p$, we obtain the \textit{algebraic bracket }$\{$ $%
,$ $\}_{k,p}:$ $J_{k}(T(M))_{p}$ $\times J_{k}(T(M))_{p}\rightarrow
J_{k-1}(T(M))_{p}.$ Note that this bracket does \textit{not} endow $%
J_{k}(T(M))_{p}$ with a Lie algebra structure. However, for $k=\infty ,$ $%
J_{\infty }(T(M))_{p}$ is a graded Lie algebra with the bracket $\{$ $,$ $%
\}_{\infty ,p},$ and is the well known Lie algebra of formal vector fields
which is extensively studied in literature ([10]). However, let $%
J_{k,0}(T(M))_{p}$ be the kernel of $\ J_{k}(T(M))_{p}\rightarrow
J_{0}(T(M))_{p}=T(M)_{p}$. Now $J_{k,0}(T(M))_{p}$ \textit{is} a Lie algebra
with the bracket $\{$ $,$ $\}_{k,p}$ which is in fact the Lie algebra of $%
\mathcal{G}_{k}(M)_{p}^{p}.$

We now define the vector bundle $J_{k}(T(M))\doteq \cup _{x\in
M}J_{k}(T(M))_{x}\rightarrow M.$ We will denote a section of $%
J_{k}(T(M))\rightarrow M$ by $\overset{(k)}{X}.$ To simplify our notation,
we will use the same notation $E$ for both the total space $E$ of a vector
bundle $E\rightarrow M$ and also for the space $\Gamma E$ of global sections
of $E\rightarrow M.$ In a coordinate system $(U,x^{i}),$ $\overset{(k)}{X}$
is of the form $\overset{(k)}{X}(x)=(x,\xi ^{i}(x),\xi _{j_{1}}^{i}(x),\xi
_{j_{2}j_{1}}^{i}(x),....,\xi _{j_{k}j_{k-1}....j_{1}}^{i}(x)),$ but we may
not have $\xi _{j_{m}j_{m-1}....j_{1}}^{i}(x)=\frac{\partial \xi
_{j_{m-1}....j_{1}}^{i}}{\partial x^{j_{m}}}(x)$ $,$ $1\leq m\leq k.$ We can
think $\overset{(k)}{X}$ also as the function $\overset{(k)}{X}(x,y)=\frac{1%
}{\alpha !}\xi _{\alpha }^{i}(x)(y-x)^{\alpha }$ where $\alpha $ is a
multi-index with $\left| \alpha \right| \leq k$ and we used summation
convention. For some $\overline{x}\in U,$ $\overset{(k)}{X}(\overline{x},y)$
is some Taylor polynomial which is \textit{not }necessarily the Taylor
polynomial of $\ \xi ^{i}(x)$ at $x=\overline{x}$ since we may \textit{not }%
have $\xi _{\alpha +j}^{i}(\overline{x})=\frac{\partial \xi _{\alpha }^{i}}{%
\partial x^{j}}(\overline{x}),$ $\left| \alpha \right| \leq k.$ Note that we
have the bundle projections $\pi _{k,m}:J_{k}(T(M))\rightarrow J_{m}(T(M))$
for $0\leq m\leq k-1,$ where $J_{0}(T(M))\doteq T(M).$ We will denote $%
J_{k}(T(M))$ by $\frak{g}_{k}(M)$ for the reason which will be clear below.

We now have the Spencer bracket $[$ $,$ $]$ defined on $\frak{g}_{k}(M)$ by
the formula

\begin{equation}
\lbrack \overset{(k)}{X},\overset{(k)}{Y}]=\{\overset{(k+1)}{X},\overset{%
(k+1)}{Y}\}+i(\overset{(0)}{X})D\overset{(k+1)}{Y}-i(\overset{(0)}{Y})D%
\overset{(k+1)}{X}\qquad k\geq 0
\end{equation}

In (2), $\overset{(k+1)}{X}$ and $\overset{(k+1)}{Y}$are arbitrary lifts of $%
\overset{(k)}{X}$ and $\overset{(k)}{Y}$, $\{$ $,$ $\}:\frak{g}%
_{k+1}(M)\times \frak{g}_{k+1}(M)\rightarrow \frak{g}_{k}(M)$ is the
algebraic bracket defined pointwise by $\{\overset{(k+1)}{X},\overset{(k+1)}{%
Y}\}(p)\doteq \{\overset{(k+1)}{X}(p),$ $\overset{(k+1)}{Y}(p)\}_{k+1,p}$
and $D$ $:\frak{g}_{k+1}(M)\rightarrow T^{\ast }\otimes \frak{g}_{k}(M)$ is
the Spencer operator given locally by the formula $(x,\xi ^{i}(x),\xi
_{j_{1}}^{i}(x),\xi _{j_{2}j_{1}}^{i}(x),....,\xi
_{j_{k+1}j_{k}....j_{1}}^{i}(x))\rightarrow (x,\frac{\partial \xi ^{i}}{%
\partial x^{j_{1}}}-\xi _{j_{1}}^{i}(x),$ $.....,\frac{\partial \xi
_{j_{k}....j_{1}}^{i}(x)}{\partial x^{j_{k+1}}}-$ $\xi
_{j_{k+1}j_{k}....j_{1}}^{i}(x)).$ Finally, the vector fields $\overset{(0)}{%
X}$ and $\overset{(0)}{Y}$ are the projections of $\overset{(k)}{X}$ and $%
\overset{(k)}{Y}$ and $i(\overset{(0)}{X})$ denotes the interior product
with respect to the vector field $\overset{(0)}{X}.$ It turns out that RHS
of (2) does not depend on the lifts $\overset{(k+1)}{X},$ $\overset{(k+1)}{Y}%
.$ The bracket $\ [$ $,$ $]$ satisfies Jacobi identity. We will refer to
[25], [26] for further details. In view of the local formulas for $\{$ $,$ $%
\}_{k,p}$ and $D$, it is elementary to make local computations using (2)
which however become formidable starting already with $k=3$. It is easy to
check that (2) gives the usual bracket formula for vector fields for $k=0.$
In fact, letting $\mathcal{X}(M)$ denote the Lie algebra of vector fields on 
$M$, we have the prolongation map $j_{k}:\mathcal{X}(M)\rightarrow \frak{g}%
_{k}(M)$ defined by $(x,\xi ^{i}(x))\rightarrow (x,\xi ^{i}(x),\partial
_{j_{1}}\xi ^{i}(x),\partial _{j_{2}j_{1}}\xi ^{i}(x),....,\partial
_{j_{k+1}j_{k}....j_{1}}\xi ^{i}(x))$ which satisfies $j_{k}\overset{(0)}{[X}%
,$ $\overset{(0)}{Y}]=[j_{k}\overset{(0)}{X},$ $j_{k}\overset{(0)}{Y}].$
Thus (2) gives the usual bracket and its derivatives when restricted to $%
j_{k}(\mathcal{X}(M)).$

Now $\frak{g}_{k}(M)$ is the transitive Lie equation in infinitesimal form $%
(TLEIF)$ determined by $\mathcal{G}_{k}(M).$ If we regard $\mathcal{G}%
_{k}(M) $ as a differentiable groupoid and construct its algebroid as in
[19], [20], we end up with $\frak{g}_{k}(M),$ justifying our notation $\frak{%
g}_{k}(M)$ for $J_{k}(T(M)).$ The projection $\pi _{k,m}:\frak{g}%
_{k}(M)\rightarrow \frak{g}_{m}(M)$ respects the bracket, i.e., it is a
homomorphism of $TLEIF^{\prime }$s.

In this way we arrive at the infinitesimal analog of (1):

\begin{equation}
......\longrightarrow \frak{g}_{k+1}(M)\longrightarrow \frak{g}%
_{k}(M)\longrightarrow ......\longrightarrow \frak{g}_{1}(M)\longrightarrow 
\frak{g}_{0}(M)\longrightarrow 0
\end{equation}

Proceeding formally, the formula $\widetilde{Diff(M)}=Lim_{k\rightarrow
\infty }\mathcal{G}_{k}(M)$ now gives $\mathcal{A}\widetilde{Diff(M)}%
=Lim_{k\rightarrow \infty }\frak{g}_{k}(M)$ where $\mathcal{A}$ stands for
the functor which assigns to a groupoid its algebroid. However note that $%
\widetilde{Diff(M)}$ is not a groupoid but rather a pseudogroup. Since a
vector field integrates to some 1-parameter group of local diffeomorphisms
(no condition on vector fields and diffeomorphisms since we have not imposed
a geometric structure yet), we naturally expect $\mathcal{A}\widetilde{%
Diff(M)}=J_{\infty }(T(M))$. As above, the vagueness in this formula is that
a vector field need not be locally determined by the Taylor expansion of its
coefficients at some point.

We now define the vector space $J_{k}(M)_{x}\doteq \{j_{k}(f)_{x}\mid f\in
C^{\infty }(M)\}$ where $C^{\infty }(M)$ denotes the set of smooth functions
on $M$ and $j_{k}(f)_{x}$ denotes the $k$-jet of $f$ at $x\in M.$ Now $%
J_{k}(M)_{x}$ is a commutative algebra over $\mathbb{R}$ with the
multiplication $\bullet $ defined by $j_{k}(f)_{x}\bullet j_{k}(g)_{x}\doteq
j_{k}(fg)_{x}.$ We define the vector bundle $J_{k}(M)\doteq $ $\cup _{x\in
M}J_{k}(M)_{x}\rightarrow M$ with the obvious differentiable structure and
projection map. The vector space of global sections of $J_{k}(M)\rightarrow
M $ is a commutative algebra with the fiberwise defined operations. We have
the projection homomorphism $\pi _{k,m}:J_{k}(M)\rightarrow J_{m}(M).$ We
will denote an element of $J_{k}(M)$ by $\overset{(k)}{f}$ which is locally
of the form $%
(x,f(x),f_{i_{1}}(x),f_{i_{2}i_{1}}(x),....,f_{i_{k}....i_{1}}(x))=(f_{%
\alpha }(x)),$ $\left| \alpha \right| \leq k.$

Now let $\overset{(k)}{X}\in \frak{g}_{k}(M)$ and $\overset{(k)}{f}\in
J_{k}(M).$ We define $\overset{(k)}{X}\overset{(k)}{f}\in J_{k}(M)$ by

\begin{equation}
\overset{(k)}{X}\overset{(k)}{f}\doteq \overset{(k)}{X}\ast \overset{(k+1)}{f%
}+i(\overset{(0)}{X})D\overset{(k+1)}{f}
\end{equation}

In (4), $\ast :$ $\frak{g}_{k}(M)\times J_{k+1}(M)\rightarrow J_{k}(M)$ is
the algebraic action of $\frak{g}_{k}(M)$ on $J_{k+1}(M)$ whose local
formula is obtained by differentiating the standard formula $\overset{(0)}{X}%
f=\xi ^{a}\partial _{a}f$ successively $k$-times and substituting jets$,$ $%
D: $ $J_{k+1}(M)\rightarrow T^{\ast }\otimes J_{k}(M)$ is the Spencer
operator defined by $%
(x,f(x),f_{i_{1}}(x),f_{i_{2}i_{1}}(x),....,f_{i_{k}....i_{1}}(x))$

$\rightarrow (x,\partial _{i_{1}}f(x)-f_{i_{1}}(x),\partial
_{i_{2}}f_{i_{1}}(x)-f_{i_{2}i_{1}}(x),....,\partial
_{i_{k}}f_{i_{k-1}....i_{1}}(x)-f_{i_{k}....i_{1}}(x))$ and $\overset{(k+1}{f%
}$ is some lift of $\overset{(k)}{f}.$ The RHS of (4) does not depend on the
lift $\overset{(k+1}{f}$. It is easy to check that $\overset{(0)}{X}\overset{%
(0)}{f}=\xi ^{a}\partial _{a}f.$ Like (2), (4) is compatible with
projection, i.e., we have $\pi _{k,m}(\overset{(k)}{X}\overset{(k)}{f})=$ $%
(\pi _{k,m}\overset{(k)}{X})($ $\pi _{k,m}\overset{(k)}{f}),$ $0\leq m\leq
k. $ Since $J_{k}(M)_{x}$ is a vector space over $\mathbb{R}$, $J_{k}(M)$ is
a module over $C^{\infty }(M).$ We will call some $\overset{(k)}{f}%
=(f_{\alpha })\in J_{k}(M)$ a smooth function if $f_{\alpha }(x)=0$ for $%
1\leq \left| \alpha \right| \leq k.$ This definition does not depend on
coordinates. Thus we have an injection $C^{\infty }(M)\rightarrow J_{k}(M)$
of algebras. If $\overset{(k)}{f}$ $\in C^{\infty }(M),$ then $\overset{(k)}{%
f}\bullet \overset{(k)}{g}=(\pi _{k,0}\overset{(k)}{f})\overset{(k)}{g}%
\doteq \overset{(0)}{f}\overset{(k)}{g}.$ Similar considerations apply to $%
J_{k}(T(M)).$

We thus obtain the $k^{th}$-order analogs of the well known formulas:

\begin{equation}
\lbrack \overset{(k)}{X},\overset{(k)}{Y}]\overset{(k)}{f}=\overset{(k)}{X}(%
\overset{(k)}{Y}\overset{(k)}{f})-\overset{(k)}{Y}(\overset{(k)}{X}\overset{%
(k)}{f})
\end{equation}

\begin{equation}
\overset{(k)}{X}(\overset{(k)}{f}\bullet \overset{(k)}{g})=(\overset{(k)}{X}%
\overset{(k)}{f})\bullet \overset{(k)}{g}+\overset{(k)}{f}\bullet (\overset{%
(k)}{X}\overset{(k)}{g})
\end{equation}

In particular, (6) gives

\begin{equation}
\overset{(k)}{X}(\overset{(0)}{f}\overset{(k)}{g})=(\overset{(0)}{X}\overset{%
(0)}{f})\overset{(k)}{g}+\overset{(0)}{f}(\overset{(k)}{X}\overset{(k)}{g})
\end{equation}

where $\overset{(0)}{X}=\pi _{k,0}\overset{(k)}{X}$. In the language of
[19], [20], (5) and (7) define a representation of the algebroid $\frak{g}%
_{k}(M)$ on the vector bundle $J_{k}(M)\rightarrow M$ (see also [25], pg.362
and [11], III, pg. 419). All constructions of this paper work also for other
such ``intrinsic'' representations of $\frak{g}_{k}(M).$ The passage from
such intrinsic representations to general representations, we believe, is
quite crucial for \textit{Erlangen Programm }and will be touched at the end
of Section 4 and in $5)$ of Section 9.

Now let $\overset{[k,r]}{\wedge }(M)_{x}$ denote the vector space of $r$%
-linear and alternating maps $\frak{g}_{k}(M)_{x}\times $ $.....\times \frak{%
g}_{k}(M)_{x}\rightarrow J_{k}(M)_{x}$ where we assume $r\geq 1,$ $k\geq 0.$
We define the vector bundle $\overset{[k,r]}{\wedge }(M)\doteq \cup _{x\in M}%
\overset{[k,r]}{\wedge }(M)_{x}\rightarrow M.$ If $\overset{[k,r]}{\omega }$ 
$\in \overset{\lbrack k,r]}{\wedge }(M)$ $(=\Gamma \overset{\lbrack k,r]}{%
\wedge }(M))$ and $\overset{(k)}{X}_{1},....,\overset{(k)}{X}_{r}\in \frak{g}%
_{k}(M),$ then $\overset{[k,r]}{\omega }(\overset{(k)}{X}_{1},....,\overset{%
(k)}{X}_{r})\in J_{k}(M)$ is defined by $(\overset{[k,r]}{\omega }(\overset{%
(k)}{X}_{1},....,\overset{(k)}{X}_{r}))(x)\doteq \overset{\lbrack k,r]}{%
\omega }(x)(\overset{(k)}{X}_{1}(x),....,$ $\overset{(k)}{X}_{r}(x))\in
J_{k}(M)_{x}.$ We define $d\overset{[k,r]}{\omega }$ by the standard
formula: $(d\overset{[k,r]}{\omega })(\overset{(k)}{X}_{1},....,\overset{(k)%
}{X}_{r+1})=$

\begin{eqnarray}
&&\frac{1}{r+1}\sum_{1\leq i\leq n+1}(-1)^{i+1}\overset{(k)}{X}_{i}\overset{%
[k,r]}{\omega }(\overset{(k)}{X}_{1},....,\overset{(i)}{\parallel },....,%
\overset{(k)}{X}_{r}) \\
&&+\frac{1}{r+1}\sum_{i\leq j+1}(-1)^{i+j}\overset{[k,r]}{\omega }([\overset{%
(k)}{X}_{i},\overset{(k)}{X}_{j}],...,\overset{(i)}{\parallel },...,\overset{%
(j)}{\parallel },...,\overset{(k)}{X}_{r+1})  \notag
\end{eqnarray}

We also define $\overset{[k,0]}{\wedge }(M)\doteq J_{k}(M)$ and $d:\overset{%
[k,0]}{\wedge }(M)\rightarrow \overset{\lbrack k,1]}{\wedge }(M)$ by $(d%
\overset{(k)}{f})\overset{(k)}{X}\doteq \overset{(k)}{X}\overset{(k)}{f}.$
We have $d:\overset{[k,r]}{\wedge }(M)\rightarrow \overset{\lbrack k,r+1]}{%
\wedge }(M)$: this follows from (6) or can be checked directly as in [9],
pg. 489, using $\overset{(k)}{f}\bullet \overset{(k)}{g}=\overset{(0)}{f}%
\overset{(k)}{g}$ if $\overset{(k)}{f}\in C^{\infty }(M).$ In view of the
Jacobi identity and the alternating character of $\overset{[k,r]}{\omega }$,
the standard computation shows $d^{2}=0.$ Thus we obtain the complex

\begin{equation}
\overset{\lbrack k,0]}{\wedge }(M)\longrightarrow \overset{\lbrack k,1]}{%
\wedge }(M)\longrightarrow \overset{\lbrack k,2]}{\wedge }(M)\longrightarrow
.....\longrightarrow \overset{\lbrack k,n]}{\wedge }(M)\qquad k\geq 0
\end{equation}

For $k=0,$ (9) gives de Rham complex.

We now assume $r\geq 1$ and define the subspace $\overset{(k,r)}{\wedge }%
(M)_{x}\subset \overset{\lbrack k,r]}{\wedge }(M)_{x}$ by the following
condition: $\overset{[k,r]}{\omega }_{x}\in \overset{\lbrack k,r]}{\wedge }%
(M)_{x}$ belongs to $\overset{(k,r)}{\wedge }(M)_{x}$ iff $\pi _{k,m}%
\overset{[k,r]}{\omega }(\overset{(k)}{X}_{1},$

$....,$ $\overset{(k)}{X}_{r})(x)\in J_{m}(M)_{x}$ depends on $\overset{(m)}{%
X}_{1}(x),....,$ $\overset{(m)}{X}_{r}(x),$ $m\leq k.$ This condition holds
vacuously for $k=0.$ Thus we obtain the projection $\pi _{k,m}:\overset{(k,r)%
}{\wedge }(M)_{x}\rightarrow \overset{(m,r)}{\wedge }(M)_{x}.$ We define the
vector bundle $\overset{(k,r)}{\wedge }(M)\doteq \cup _{x\in M}\overset{(k,r)%
}{\wedge }(M)_{x}\rightarrow M$ and set $\overset{(k,0)}{\wedge }=\overset{%
[k,0]}{\wedge }=J_{k}(M).$

\begin{definition}
An exterior $(k,r)$-form $\overset{(k,r)}{\omega }$ on $M$ is a smooth
section of the vector bundle $\overset{(k,r)}{\wedge }(M)\rightarrow M.$
\end{definition}

An explicit description of $\overset{(k,r)}{\omega }$ in local coordinates
is not without interest but will not be done here. Applying $\pi _{k,m}$ to
(8), we deduce $d:\overset{(k,r)}{\wedge }(M)\rightarrow \overset{(k,r+1)}{%
\wedge }(M)$ and the commutative diagram

\begin{equation}
\begin{array}{ccc}
\overset{(k+1,r)}{\wedge } & \overset{d}{\longrightarrow } & \overset{%
(k+1,r+1)}{\wedge } \\ 
\downarrow _{\pi } &  & \downarrow _{\pi } \\ 
\overset{(k,r)}{\wedge } & \overset{d}{\longrightarrow } & \overset{(k,r+1)}{%
\wedge }
\end{array}
\end{equation}

Thus we obtain the array

\begin{equation}
\begin{array}{ccccccccc}
\overset{(\infty ,0)}{\wedge } & \longrightarrow & \overset{(\infty ,1)}{%
\wedge } & \longrightarrow & \overset{(\infty ,2)}{\wedge } & \longrightarrow
& ..... & \longrightarrow & \overset{(\infty ,n)}{\wedge } \\ 
..... &  & ..... &  & ..... &  & ..... &  & ..... \\ 
\downarrow &  & \downarrow &  & \downarrow &  & \downarrow &  & \downarrow
\\ 
\overset{(2,0)}{\wedge } & \longrightarrow & \overset{(2,1)}{\wedge } & 
\longrightarrow & \overset{(2,2)}{\wedge } & \longrightarrow & ..... & 
\longrightarrow & \overset{(2,n)}{\wedge } \\ 
\downarrow &  & \downarrow &  & \downarrow &  & \downarrow &  & \downarrow
\\ 
\overset{(1,0)}{\wedge } & \longrightarrow & \overset{(1,1)}{\wedge } & 
\longrightarrow & \overset{(1,2)}{\wedge } & \longrightarrow & ..... & 
\longrightarrow & \overset{(1,n)}{\wedge } \\ 
\downarrow &  & \downarrow &  & \downarrow &  & \downarrow &  & \downarrow
\\ 
\overset{(0,0)}{\wedge } & \longrightarrow & \overset{(0,1)}{\wedge } & 
\longrightarrow & \overset{(0,2)}{\wedge } & \longrightarrow & ..... & 
\longrightarrow & \overset{(0,n)}{\wedge }
\end{array}
\end{equation}

where all horizontal maps are given by $d$ and all vertical maps are
projections. The top sequence is defined algebraically by taking projective
limits in each coloumn.

Using $\bullet $, we also define the wedge product of $(k,r)$-forms in each
row by the standard formula which turns $\overset{(k,\ast )}{\wedge }%
(M)\doteq \oplus _{0\leq r\leq n}\overset{(k,r)}{\wedge }(M)$ into an
algebra. This algebra structure descends to the cohomology. In particular,
we obtain an algebra structure on the cohomology of the top row of (11),
which we will denote by $H^{\ast }(\frak{g}_{\infty }(M),J_{\infty }(M))$.

Now let $C^{k+1,r}(M)\doteq $ $Kernel$ $(\pi _{\infty ,k}:\overset{(\infty
,r)}{\wedge }(M)\rightarrow \overset{(k,r)}{\wedge }(M))$ and $\mathcal{C}%
^{k+1,\ast }(M)\doteq \oplus _{0\leq r\leq n}C^{k+1,r}(M),$ $k\geq 0.$ This
gives the array

\begin{equation}
\begin{array}{ccccccccc}
\overset{(\infty ,0)}{\wedge }(M) & \longrightarrow & \overset{(\infty ,1)}{%
\wedge }(M) & \longrightarrow & \overset{(\infty ,2)}{\wedge }(M) & 
\longrightarrow & ..... & \longrightarrow & \overset{(\infty ,n)}{\wedge }(M)
\\ 
\uparrow &  & \uparrow &  & \uparrow &  & \uparrow &  & \uparrow \\ 
C^{1,0}(M) & \longrightarrow & C^{1,1}(M) & \longrightarrow & C^{1,2}(M) & 
\longrightarrow & ..... & \longrightarrow & C^{1,n}(M) \\ 
\uparrow &  & \uparrow &  & \uparrow &  & \uparrow &  & \uparrow \\ 
C^{2,0}(M) & \longrightarrow & C^{2,1}(M) & \longrightarrow & C^{2,2}(M) & 
\longrightarrow & .... & \longrightarrow & C^{2,n}(M) \\ 
\uparrow &  & \uparrow &  & \uparrow &  & \uparrow &  & \uparrow \\ 
..... &  & ..... &  & ..... &  & ..... &  & ..... \\ 
C^{\infty ,0}(M) & \longrightarrow & C^{\infty ,1}(M) & \longrightarrow & 
C^{\infty ,2}(M) & \longrightarrow & ..... & \longrightarrow & C^{\infty
,n}(M)
\end{array}
\end{equation}

where all horizontal maps in (12) are restrictions of $d$ and all vertical
maps are inclusions. The filtration in (12) preserves wedge product. In
fact, we have $C^{k,r}(M)\wedge $ $C^{s,t}(M)\subset C^{k+s,r+t}(M)$ which
follows easily from the definition of $\bullet .$ We will denote the
spectral sequence of algebras determined by the filtration in (12) by $%
\mathcal{U}_{M}$ and call $\mathcal{U}_{M}$ the universal spectral sequence
of $M.$

The above construction will be relevant in the next section. However, we
remark here that (11) contains no information other than $H_{deR}^{\ast }(M,%
\mathbb{R)}.$ To see this, we observe that if $\overset{(k)}{X}\overset{(k)}{%
f}=0$ for all $\overset{(k)}{X}\in \frak{g}_{k}(M))$, then $\overset{(k)}{f}%
\in \mathbb{R}\subset C^{\infty }(M).$ Thus the kernel of the first
operators in the horizontal rows of (11) define the constant sheaf $\mathbb{R%
}$ on $M.$ Since $\overset{(k,r)}{\wedge }(M)$ is a module over $C^{\infty
}(M)$, each row of (10) (which is easily shown to be locally exact) is a
soft resolution of the constant sheaf $\mathbb{R}$ and thus computes $%
H_{deR}^{\ast }(M,\mathbb{R)}.$

\section{Homogeneous geometries}

\begin{definition}
A homogeneous geometry on a differentiable manifold $M$ is a diagram

\begin{equation}
\begin{array}{ccccccccc}
..... & \longrightarrow  & \mathcal{G}_{2}(M) & \longrightarrow  & \mathcal{G%
}_{1}(M) & \longrightarrow  & M\times M & \longrightarrow  & 1 \\ 
\uparrow  &  & \uparrow  &  & \uparrow  &  & \parallel  &  &  \\ 
..... & \longrightarrow  & \mathcal{S}_{2}(M) & \longrightarrow  & \mathcal{S%
}_{1}(M) & \longrightarrow  & M\times M & \longrightarrow  & 1
\end{array}
\end{equation}

where $i)$ $\mathcal{S}_{k}(M)$ is a $TLEFF$ for all $k\geq 1$ and therefore 
$\mathcal{S}_{k}(M)$ $\subset \mathcal{G}_{k}(M)$ and the vertical maps are
inclusions $ii)$ The horizontal maps in the bottom sequence of (13) are
restrictions of the projection maps in the top sequence and are surjective
morphisms.
\end{definition}

With an abuse of notation, we will denote (13) by $\mathcal{S}_{\infty }(M)$
and call $\mathcal{S}_{\infty }(M)$ a homogeneous geometry on $M.$ We thus
imagine that the lower sequence of (13) ``converges'' to some (pseudo)group $%
G\subset \widetilde{Diff(M)}$ which acts transitively on $M.$ However, $G$
may be far from being a Lie group and it may be intractible to deal with $G$
directly. The idea of Definition 2 is to work with the arrows of $G$ rather
than to work with $G$ itself.

\begin{definition}
Let $\mathcal{S}_{\infty }(M)$ be a homogeneous geometry on $M.$ $\mathcal{S}%
_{\infty }(M)$ \ is called a Klein geometry if there exists an integer $%
m\geqslant 1$ with the following property: If $(f_{m})_{y}^{x}\in \mathcal{S}%
_{m}(M)_{y}^{x}$, then there exists a unique local diffeomorphism $g$ with $%
g(x)=y$ satisfying $i)$ $j_{k}(g)_{y}^{x}=(f_{m})_{y}^{x}$ \ $ii)$ $%
j_{k}(g)_{g(z)}^{z}\in \mathcal{S}_{m}(M)_{g(z)}^{z}$ for all $z$ near $x.$
The smallest such integer (uniquely determined if $M$ is connected) is
called the order of the Klein geometry.
\end{definition}

In short, a Klein geometry is a transitive pseudogroup whose local
diffeomorphisms are uniquely determined by any of their $m$-arrows and we
require $m$ to be the smallest such integer. Once we have a Klein geometry $%
\mathcal{S}_{\infty }(M)$ of order $m,$ then all $\mathcal{S}_{k}(M),$ $%
k\geqslant m+1$ are uniquely determined by $\mathcal{S}_{m}(M)$ and $%
\mathcal{S}_{k+1}(M)\rightarrow \mathcal{S}_{k}(M)$ is an isomorphism for
all $k\geqslant m,$ i.e., the Klein geometry $\mathcal{S}_{\infty }(M)$
prolongs in a unique way to a homogeneous geometry and all the information
is contained in terms up to order $m.$ We will thus denote a Klein geometry
by $\mathcal{S}_{(m)}(M).$ We will take a closer look at these geometries in
the next section. For instance, let $M^{2n}$ be a complex manifold. We
define $\mathcal{S}_{k}(M)_{y}^{x}$ by the following condition: $%
(f_{k})_{y}^{x}\in \mathcal{G}_{k}(M)_{y}^{x}$ belongs to $\mathcal{S}%
_{k}(M)_{y}^{x}$ if there exists a local holomorphic diffeomorphism $g$ with 
$g(x)=y$ and $j^{k}(g)_{y}^{x}=(f_{k})_{y}^{x}.$ We see that $TLEFF$ $%
\mathcal{S}_{k}(M)$ defined by $\mathcal{S}_{k}(M)\doteq \cup _{x,y\in M}%
\mathcal{S}_{k}(M)_{y}^{x}$ satisfies conditions of Definition 3 and
therefore a complex structure determines a homogeneous geometry which is not
necessarily a Klein geometry. Similarly, a symplectic or contact structure
determines a homogeneous geometry (since these structures have no local
invariants) which need not be Klein. More generally, any transitive
pseudogroup determines a homogeneous geometry via its arrows.

Now given a homogeneous geometry $\mathcal{S}_{\infty }(M)$, we will sketch
the construction of its infinitesimal geometry $\frak{s}_{\infty }(M)$
referring to [25], [26] for further details. Let $x\in M$ and $X$ be a
vector field defined near $x.$ Let $f_{t}$ be the 1-parameter group of local
diffeomorphisms generated by $X.$ Suppose $X$ has the property that $%
j_{k}(f_{t})_{y_{t}}^{x}$ belongs to $\mathcal{S}_{k}(M)_{y_{t}}^{x}$ for
all small $t$ with $t\geqslant 0$ where $y_{t}=f_{t}(x)$. This is actually a
condition only on the $k$-jet of $X$ at $x.$ In this way we define the
subspace $\frak{s}_{k}(M)_{x}\subset \frak{g}_{k}(M)_{x}$ which consists of
those $\overset{(k)}{X}(x)$ satisfying this condition. We define the vector
subbundle $\frak{s}_{k}(M)\doteq \cup _{x\in M}\frak{s}_{k}(M)_{x}%
\rightarrow M$ of $\frak{g}_{k}(M)\rightarrow M$ and the bracket (2) on $%
\frak{g}_{k}(M)$ restricts to a bracket on $\frak{s}_{k}(M)$ and $\frak{s}%
_{k}(M)$ is the $TLEIF$ determined by $\mathcal{S}_{k}(M).$

In this way we arrive at the diagram

\begin{equation}
\begin{array}{ccccccccccc}
.... & \longrightarrow & \frak{g}_{k}(M) & \longrightarrow & .... & 
\longrightarrow & \frak{g}_{1}(M) & \longrightarrow & \frak{g}_{0}(M) & 
\longrightarrow & 0 \\ 
&  & \uparrow &  & \uparrow &  & \uparrow &  & \uparrow &  &  \\ 
.... & \longrightarrow & \frak{s}_{k}(M) & \longrightarrow & .... & 
\longrightarrow & \frak{s}_{1}(M) & \longrightarrow & \frak{s}_{0}(M) & 
\longrightarrow & 0
\end{array}
\end{equation}

where the bottom horizontal maps are restrictions of the upper horizontal
maps and are surjective morphisms. The vertical maps are injective morphisms
induced by inclusion. Thus (14) is the infinitesimal analog of (13). We will
denote (14) by $\frak{s}_{\infty }(M)$ and call $\frak{s}_{\infty }(M)$ the
infinitesimal geometry of $\mathcal{S}_{\infty }(M).$

Now we define $L_{\overset{(k)}{Y}}:\overset{(k,r)}{\wedge }\rightarrow 
\overset{(k,r)}{\wedge }$ and $i_{\overset{(k)}{Y}}:\overset{(k,r+1)}{\wedge 
}\rightarrow \overset{(k,r)}{\wedge }$ by the standard formulas: $(L_{%
\overset{(k)}{Y}}\overset{(k,r)}{\omega })(\overset{(k)}{X}_{1},....,%
\overset{(k)}{X}_{r})\doteq \overset{(k)}{Y}(\overset{(k,r)}{\omega }(%
\overset{(k)}{X}_{1},....,\overset{(k)}{X}_{r}))-\overset{(k,r)}{\omega }([%
\overset{(k)}{Y},\overset{(k)}{X}_{1}],\overset{(k)}{X}_{2}....,$

$\overset{(k)}{X}_{r})-....-\overset{(k,r)}{\omega }([\overset{(k)}{X}_{1},%
\overset{(k)}{X}_{2}....,[\overset{(k)}{Y},\overset{(k)}{X}_{r}])$ and $(i_{%
\overset{(k)}{Y}}\overset{(k+1,r)}{\omega })(\overset{(k)}{X}_{1},....,%
\overset{(k)}{X}_{r})\doteq (r+1)$ $\ \overset{(k+1,r)}{\omega }(\overset{(k)%
}{Y},\overset{(k)}{X}_{1},....,\overset{(k)}{X}_{r}).$ We also define $(\pi
_{k,m}L_{\overset{(k)}{Y}}):\overset{(m,r)}{\wedge }\rightarrow \overset{%
(m,r)}{\wedge }$ by $\pi _{k,m}L_{\overset{(k)}{Y}}\doteq L_{\pi _{k,m}%
\overset{(k)}{Y}}$ and similarly $\pi _{k,m}i_{\overset{(k)}{Y}}=i_{\pi
_{k,m}\overset{(k)}{Y}}.$ With these definitions we obtain the well known
formulas

\begin{equation}
L_{\overset{(k)}{X}}=d\circ i_{\overset{(k)}{X}}+i_{\overset{(k)}{X}}\circ d
\end{equation}

\begin{equation}
L_{\overset{(k)}{X}}\circ d=d\circ L_{\overset{(k)}{X}}
\end{equation}

We will now indicate briefly how a homogeneous geometry $\mathcal{S}_{\infty
}(M)$ gives rise to various spectral sequences.

$1)$ We define $\overset{(k,r)}{\frak{s}}(M)\subset \overset{(k,r)}{\wedge }%
(M)$ by the following condition: Some $\overset{(k,r)}{\omega }$ belongs to $%
\overset{(k,r)}{\frak{s}}(M)$ if $\quad i)$ $L_{\overset{(k)}{X}}\overset{%
(k,r)}{\omega }=0,$ $\overset{(k)}{X}\in \frak{s}_{k}(M)$ $\quad ii)$ $i_{%
\overset{(k)}{X}}\overset{(k,r)}{\omega }=0,$ $\overset{(k)}{X}\in \frak{s}%
_{k}(M).$ (15) and (16) show that $d:\overset{(k,r)}{\frak{s}}(M)\rightarrow 
\overset{(k,r+1)}{\frak{s}}(M)$ and (10) holds. Thus we arrive at (11) and
(12). Recall that if $\frak{g}$ is any Lie algebra with a representation on $%
V$ and $\frak{s}\subset \frak{g}$ is a subalgebra, then we can define the
relative Lie algebra cohomology groups $H^{\ast }(\frak{g},\frak{s},V)$ ([6])%
$.$ Since our construction is modelled on the definition of $H^{\ast }(\frak{%
g},\frak{s},V),$ we will denote the cohomology of the top row of (11) by $%
H^{\ast }(\frak{g}_{\infty }(M),\frak{s}_{\infty }(M),J_{\infty }(M))$ in
this case.

$2)$ We will first make

\begin{definition}
$\Theta (\mathcal{S}_{k}(M))\doteq \{\overset{(k)}{f}\in J_{k}(M)\mid 
\overset{(k)}{X}\overset{(k)}{f}=0$ for all $\overset{(k)}{X}\in \frak{s}%
_{k}(M)\}$
\end{definition}

(6) shows that $\Theta (\mathcal{S}_{k}(M))$ is a subalgebra of $J_{k}(M).$
\ We will call $\Theta (\mathcal{S}_{k}(M))$ the $k^{th}$ order structure
algebra of the homogeneous geometry $\mathcal{S}_{\infty }(M).$ Note that $%
\Theta (\mathcal{G}_{k}(M))=\mathbb{R}.$ For $k\geq 1,$ we define $\overset{%
[k,r]}{\frak{s}}(M)_{x}$ as the space of alternating maps $\frak{s}%
_{k}(M)_{x}\times ....\times \frak{s}_{k}(M)_{x}\rightarrow J_{k}(M)_{x}$
and define $\overset{[k,r]}{\frak{s}}(M)$ as in Section 2. We have the
restriction map $\theta _{(k,r)}:\overset{[k,r]}{\wedge }(M)\rightarrow 
\overset{\lbrack k,r]}{\frak{s}}(M)$ whose kernel will be denoted by $%
\overset{[k,r]}{\wedge }_{\frak{s}}(M).$ Since $\Theta (\mathcal{S}%
_{k}(M))=Ker(\theta _{(k,1)}\circ d),$ we obtain the commutative diagram

\begin{equation}
\begin{array}{ccccccccccc}
&  & 0 &  & 0 &  & 0 &  & .... &  & 0 \\ 
&  & \downarrow &  & \downarrow &  & \downarrow &  & \downarrow &  & 
\downarrow \\ 
&  & \Theta (\mathcal{S}_{k}(M)) & \longrightarrow & \overset{[k,1]}{\wedge }%
_{\frak{s}}(M) & \longrightarrow & \overset{[k,2]}{\wedge }_{\frak{s}}(M) & 
\longrightarrow & .... & \longrightarrow & \overset{[k,n]}{\wedge }_{\frak{s}%
}(M) \\ 
&  & \downarrow &  & \downarrow &  & \downarrow &  & \downarrow &  & 
\downarrow \\ 
\mathbb{R} & \longrightarrow & \overset{[k,0]}{\wedge }(M) & \longrightarrow
& \overset{[k,1]}{\wedge }(M) & \longrightarrow & \overset{[k,2]}{\wedge }(M)
& \longrightarrow & .... & \longrightarrow & \overset{[k,2]}{\wedge }(M) \\ 
&  & \downarrow &  & \downarrow &  & \downarrow &  & \downarrow &  & 
\downarrow \\ 
0 & \longrightarrow & \frac{\overset{[k,0]}{\wedge }(M)}{\Theta (\mathcal{S}%
_{k}(M))} & \longrightarrow & \overset{[k,1]}{\frak{s}}(M) & \longrightarrow
& \overset{[k,2]}{\frak{s}}(M) & \longrightarrow & .... & \longrightarrow & 
\overset{[k,n]}{\frak{s}}(M) \\ 
&  & \downarrow &  & \downarrow &  & \downarrow &  & \downarrow &  & 
\downarrow \\ 
&  & 0 &  & 0 &  & 0 &  & .... &  & 0
\end{array}
\end{equation}

We will call (17) the horizontal crossection of the representation triple $(%
\mathcal{G}_{\infty }(M),$

$\mathcal{S}_{\infty }(M),J_{\infty }(M))$ at order $k$. Passing to the long
exact sequence in (17), we see that local cohomology of the top and bottom
sequences coincide (with a shift in order) in view of the local exactness of
the middle row. Now defining $\overset{(k,r)}{\frak{s}}(M)$ as in Section 2,
(17) defines an exact sequence of three spectral sequences where the middle
spectral sequence is (12). Note that local exactness of the top and bottom
sequences would imply that their cohomologies coincide with the sheaf
cohomology groups $H^{\ast }(M,\Theta (\mathcal{S}_{k}(M))$ and $H^{\ast }(M,%
\frac{\overset{[k,0]}{\wedge }(M)}{\Theta (\mathcal{S}_{k}(M))})$
respectively since partition of unity applies to the spaces in these
sequences. We will denote the limiting cohomology of the top and bottom
sequences in (17) respectively by $H^{\ast }(\frak{s}_{\infty }(M),0)$ and $%
H^{\ast }(\frak{s}_{\infty }(M),J_{\infty }(M)).$ The reader may compare
(17) to the diagram on page 183 in [25] which relates Spencer sequence to
Janet sequence, called $\mathcal{P}$-sequence in [25].

Before we end this section, note that $(f_{k+1})_{y}^{x}\in \mathcal{S}%
_{k}(M)_{y}^{x}$ defines an isomorphism $(f_{k+1})_{y}^{x}:J_{k}(M)_{x}%
\rightarrow J_{k}(M)_{y}\ $(in fact, $(f_{k})_{y}^{x}$ does it) and also an
isomorphism $(f_{k+1})_{y}^{x}:\frak{s}_{k}(M)_{x}\rightarrow \frak{s}%
_{k}(M)_{y}.$ Let us assume that $\mathcal{S}_{\infty }(M)$ is defined by
some $G\subset Diff(M)$ as in the case of a symplectic structure,
homogeneous complex structure or a Klein geometry (see Section 4). Now $G$
acts on $\overset{(k,r)}{\frak{s}}(M)$ (defined as in 1) or 2) above) by $(g%
\overset{(k,r)}{\omega })(\overset{(k)}{X}_{1},....,\overset{(k)}{X}%
_{r})(p)\doteq j_{k+1}(g)_{p}^{q}(\overset{(k,r)}{\omega }%
(q)(j_{k+1}(g^{-1})_{q}^{p}\overset{(k)}{X}%
_{1}(p),....,j_{k+1}(g^{-1})_{q}^{p}\overset{(k)}{X}_{r}(p))$ where $g(q)=p.$
The cochains which are invariant under this action form a subcomplex whose
cohomology can be localized, i.e., can be computed at any point of $M.$ If $%
\mathcal{S}_{\infty }(M)=\mathcal{G}_{\infty }(M)$ and $G=Diff(M)$, then an
invariant form must vanish but this need not be the case for a homogeneous
geometry. We will not go into the precise description of this cohomology
here though it is quite relevant for Klein geometries in Section 4 and can
be expressed in terms of some relative Lie algebra cohomology groups in this
case.

\section{Klein geometries}

Let $G$ be a Lie group and $H$ a closed subgroup. $G$ acts on the left coset
space $G/H$ on the left. Let $o$ denote the coset of $H.$ Now $H$ fixes $o$
and therefore acts on the tangent space $T(G/H)_{o}$ at $o$. However some
elements of $H$ may act as identity on $T(G/H)_{o}.$ The action of $h\in H$
(we regard $h$ as a transformation and use the same notation) on $T(G/H)_{o}$
depends only on $1$-arrow of $h$ with source and target at $o.$ Let $%
H_{1}\subset H$ be the normal subgroup of $H$ consisting of elements which
act as identity on $T(G/H)_{o}.$ To recover $H_{1},$ we consider $1$-jets of
vector fields at $o$ which we will denote by $J_{1}T(G(H)_{o}.$ The action
of $h$ on $J_{1}T(G(H)_{o}$ depends only on $2$-arrow of $h.$ Now some
elements $h\in H_{1}$ may still act as identity at $J_{1}T(G(H)_{o}$ and we
define the normal subgroup $H_{2}\subset $ $H_{1}$ consisting of those
elements. Iterating this procedure, we obtain a decreasing sequence of
normal subgroups $\{1\}\subset ...\subset H_{k}\subset H_{k-1}\subset
.....\subset H_{2}\subset H_{1}\subset H_{0}=H$ which stabilizes at some
group $N$ which is the largest normal subgroup of $G$ contained in $H$ $($%
see [29], pg. 161). We will call the smallest integer $m$ satisfying $%
N=H_{m} $ the order of the Klein pair $(G,H)$. In this case, it is easy to
show that $g\in G$ is uniquely determined modulo $N$ by any of its $m$%
-arrows. $G$ acts effectively on $G/H$ iff $N=\{1\}.$ If $(G,H)$ is a Klein
pair of order $m,$ then so is $(G/N,H/N)$ which is further effective. We
will call $N$ the ghost of the Klein pair $(G,H)$ since it can not be
detected from the action and therefore \textit{may} have implications that
fall outside the scope of \textit{Erlangen Programm. }We will touch this
issue again in $5)$ of Section 9. Thus we see that an \textit{effective }%
Klein pair $(G,H)$ of order $m$ determines a Klein geometry $\mathcal{S}%
_{(m)}(G/H)$ according to Definition 3 where the local diffeomorphisms
required by Definition 3 are restrictions of global diffeomorphisms of $G/H$
which are induced by the elements of $G.$ Conversely, let $\mathcal{S}%
_{(m)}(M)$ be a Klein geometry according to Definition 3 and let $\widetilde{%
M}$ be the universal covering space of $M.$ We pull back the pseudogroup on $%
M$ to a pseudogroup on $\widetilde{M}$ using the local diffeomorphism $\pi :%
\widetilde{M}\rightarrow M.$ Using simple connectedness of $\widetilde{M}$
and a mild technical assumption which guarantees that the domains of the
local diffeomorphisms do not become arbitrarily small, the standard
monodromy argument shows that a local diffeomorphism defined on $\widetilde{M%
}$ in this way uniquely extends to some global diffeomorphism on $\widetilde{%
M}.$ This construction is essentially the same as the one given in [30] on
page 139-146. The global diffeomorphisms on $\widetilde{M}$ obtained in this
way form a Lie group $G.$ If $H\subset G$ is the stabilizer of some $p\in 
\widetilde{M},$ then $H$ is isomorphic to $\mathcal{S}_{(m)}(M)_{q}^{q}$
where $\pi (p)=q.$ To summarize, a Klein geometry $\mathcal{S}_{(m)}(M)$
according to Definition 3 becomes an effective Klein pair $(G,H)$ of order $%
m $ when pulled back to $\widetilde{M}$. Conversely, an effective Klein pair 
$(G,H)$ of order $m$ defines a Klein geometry $\mathcal{S}_{(m)}(M)$ if we
mode out by the action of a discrete subgroup. Keeping this relation in
mind, we will consider an effective Klein pair $(G,H)$ as our main object
below.

Now the above filtration of normal subgroups gives the diagram

\begin{equation}
\begin{array}{ccccccccccc}
\mathcal{G}_{m}(G/H)_{o}^{o} & \longrightarrow  & \mathcal{G}%
_{m-1}(G/H)_{o}^{o} & \longrightarrow  & .... & \longrightarrow  & \mathcal{G%
}_{1}(G/H)_{o}^{o} & \longrightarrow  & 1 &  &  \\ 
\uparrow  &  & \uparrow  &  & \uparrow  &  & \uparrow  &  &  &  &  \\ 
H & \longrightarrow  & H/H_{m} & \longrightarrow  & ... & \longrightarrow  & 
H/H_{1} & \longrightarrow  & 1 &  & 
\end{array}
\end{equation}

where the spaces in the top sequence are jet groups in our universal
situation in Section 2 and the vertical maps are injections. Since the
kernels in the upper sequence are vector groups, we see that $H_{1}$ is
solvable. As before, we now define $\mathcal{S}_{m}(M)_{y}^{x}$ which
consists of $m$-arrows of elements of $G$ and define $\mathcal{S}%
_{m}(M)\doteq \cup _{x,y\in M}\mathcal{S}_{m}(M)_{y}^{x}.$ As in the case $%
Diff(M)\times M\rightarrow M$ in Section 2, we obtain the map $G\times
G/H\rightarrow \mathcal{S}_{m}(M)$ defined by $(g,x)\rightarrow $ $m$-arrow
of $g$ starting at $x$ and ending at $g(x).$ This map is surjective by
definition and this time also injective by the definition of $m.$ Thus we
obtain a concrete realization of $\mathcal{S}_{m}(M)$ as $G\times G/H.$ Note
that $\mathcal{S}_{m}(M)_{o}^{o}=H.$ Going downwards in the filtration, we
obtain the commutative diagram

\begin{equation}
\begin{array}{ccc}
G\times G/H & \longrightarrow & \mathcal{S}_{m}(M) \\ 
\downarrow &  & \downarrow \\ 
G/H_{m-1}\times G/H & \longrightarrow & \mathcal{S}_{m-1}(M) \\ 
\downarrow &  & \downarrow \\ 
.... & \longrightarrow & ... \\ 
\downarrow &  & \downarrow \\ 
G/H_{2}\times G/H & \longrightarrow & \mathcal{S}_{2}(M) \\ 
\downarrow &  & \downarrow \\ 
G/H_{1}\times G/H & \longrightarrow & \mathcal{S}_{1}(M)
\end{array}
\end{equation}

For instance, the bottom map in (19) is defined by $\{xH_{1}\}\times
\{yH\}\rightarrow 1$-arrow of the diffeomorphism $\{xH_{1}\}$ starting at
the coset $\{yH\}$ and ending at the coset $\{xyH\}.$ Note that this is not
a group action since $G/H_{1}$ is not a group but the composition and
inversion of $1$-arrows are well defined. This map is a bijection by the
definition of $H_{1}.$ Fixing one such $1$-arrow$,$ all other $1$-arrows
starting at $\{yH\}$ are generated by composing this arrow with elements of $%
\mathcal{S}_{1}(M)_{y}^{y}=I_{xy^{-1}}H/I_{xy^{-1}}H_{1}$ where $I_{xy^{-1}}$
is the inner automorphism of $G$ determined by $xy^{-1}\in G.$ The vertical
projections on the right coloumn of (19) are induced by projection of jets
as in Sections 2, 3 and the projections on the left coloumn are induced by
projections on the first factor and identity map on the second factor.

A Lie group $G$ is clearly an effective Klein pair $(G,\{1\})$ with order $%
m=0.$ For many Klein geometries we have $m=1$. This is the case, for
instance, if $H$ is compact, in which case we have an invariant metric, $H$
is discrete or $(G,H)$ is a reductive pair which is extensively studied in
literature from the point of view of principal bundles ([12]). If $G$ is
semisimple, it is not difficult to show that the order of $(G,H)$ is at most
two (see [18], pg. 131). For instance, let $M$ be a homogeneous complex
manifold, i.e., $Aut(M)$ acts transitively on $M.$ If $M$ is compact, then $%
M=G/H$ for some complex Lie group $G$ and a closed complex subgroup $H$
([34]). If futher $\pi _{1}(M)$ is finite, then $G/H=\overline{G}/\overline{H%
}$ as complex manifolds for some semisimple Lie group $\overline{G}$ ([34]).
Thus it follows that jets of order greater than two do not play any role in
the complex structure of $M$ in this case. If $G$ is reductive, it is stated
in [35] that the order of $(G,H)$ is at most three. On the other hand, for
any positive integer $m$, an effective Klein pair $(G;H)$ of order $m$ is
constructed in [1] such that $G/H$ is open and dense in some weighted
projective space. Other examples of Klein pairs of arbitrary order are
communicated to us by the author of [35]. However, we do not know the answer
to the following question

$\mathbf{Q1:}$ For some positive integer $m,$ does there exist a Klein pair $%
(G,H)$ of order $m$ such that $G/H$ compact?

It is crucial to observe $\ $that $(G_{1},H_{1})$ and $(G_{2},H_{2})$ may be
two Klein pairs with different orders with $G_{1}/H_{1}$ homeomorphic to $%
G_{2}/H_{2}.$ For instance, let $H\subset G$ be complex Lie groups with $\pi
_{0}(G)=\pi _{1}(G/H)=0$ and $G/H$ compact. If $M$ $\subset G$ is a maximal
compact subgroup, then $(M,M\cap H)$ is a Klein pair of order one and $%
G/H=M/M\cap H$ as topological manifolds (in fact, $G/H$ is Kaehler iff its
Euler characteristic is nonzero, see [34], [4]). The crucial fact here is
that an abundance of Lie groups may act transitively on a homogeneous space $%
M$ with different orders and the topology (but not the analytic structure)
of $M$ is determined by actions of order one only (the knowledge of a
particular such action suffices, see Theorem VI in [34] where a detailed
description of complex homogeneous spaces is given. It turns out that ``they
are many more than we expect'' as stated there).

Now let $(G,H)$ be an effective Klein pair of order $m$ and let $\mathcal{L}%
(G)$ be the Lie algebra of $G.$ We have the map $\sigma :\mathcal{L}%
(G)\rightarrow J_{m}(T(G/H))_{p}$ defined by $X\rightarrow j_{m}(X^{\ast
})_{p}$ where $X^{\ast }$ is the vector field on $G/H$ induced by $X$ and $%
p\in G/H.$ $\sigma $ is a homomorphism of Lie algebras where the bracket on $%
J_{m}(T(G/H))_{o}$ is the algebraic bracket defined in Section 2. Note that
the map $X\rightarrow X^{\ast }$ is injective due to effectiveness and also
surjective due to transitivity. It is now easy to give a description of the
infinitesimal analog of (19). We can thus express everything defined in
Sections 2, 3 in concrete terms which will enable us to use the highly
developed structure theory of (semisimple) Lie groups ([17]). A detailed
description of $\frak{s}_{m}(M)$ is given in [26] (Theorem 15 on pg. 199).
The formula on pg. 200 in [26] is the same as the formula (15), Example
3.3.7, pg. 104 in the recent book [20], but higher order jets and Spencer
operator remain hidden in (15) and also in (4), Example 3.2.9, pg. 98 in
[20].

The following situation deserves special mention: Let $G$ be a complex Lie
group, $H$ a closed complex subgroup with a holomorphic representation on
the vector space $V$ and let $E\rightarrow G/H$ be the associated
homogeneous vector bundle of $G\rightarrow G/H.$ We now have the sheaf
cohomology groups $H^{\ast }(G/H,\mathbf{E})$ where $\mathbf{E}$ denotes the
sheaf of holomorphic sections of $E\rightarrow G/H.$ Borel-Weil theorem is
derived in [4] from $H^{0}(G/H,\mathbf{E)}$. If the Klein pair $(G,H)$ has
order $m$ and is effective, the principal bundle $G\rightarrow G/H$ can be
identified with the principal bundle $\mathcal{S}_{m}(M)^{(o)}\rightarrow
G/H $ where $\mathcal{S}_{m}(M)^{(o)}$ consists of $m$-arrows in $G/H$ with
source at the coset $o$ of $H$ (see Section 7). If $gx=y$, $g\in G,$ $x,y\in
G/H$, then the $m$-arrow of $g$ gives an isomorphism $E_{x}\rightarrow E_{y}$
between the fibers. Consequently, the action of $G$ on sections of $%
E\rightarrow G/H$ is equivalent to the representation of the $TLEFF$ $%
\mathcal{S}_{m}(M)=G\times G/H$ on $E\rightarrow G/H.$ On the infinitesimal
level, this gives a representation of $\frak{s}_{m}(M)$ on $E\rightarrow G/H$
and we can also define the cohomology groups $H^{\ast }(\frak{s}_{m}(M),E)$
as in [19], [20]. Letting $\Theta $ denote the sections of $E$ killed by $%
\frak{s}_{m}(M),$ we see that these two cohomology groups are related by a
diagram similar to (17). It is crucial to observe here how the order of jets
remains hidden in $H^{\ast }(G/H,\mathbf{E})$.

\section{Truncated geometries}

\begin{definition}
Some $TLEFF$ $\mathcal{S}_{k}(M)$ is called a truncated geometry on $M$ of
order $k.$
\end{definition}

We will view a truncated geometry $\mathcal{S}_{k}(M)$ as a diagram (13)
where all $\mathcal{S}_{m}(M)$ for $m\leq k$ are defined as projections of $%
\mathcal{S}_{k}(M)$ and $\mathcal{S}_{m}(M)$ for $m\geqslant k+1$ do not
exist. A homogeneous geometry defines a truncated geometry of any order. The
question arises whether some truncated geometry always prolongs uniquely to
some homogeneous geometry. The answer turns out to be negative. For
instance, let $(M,g)$ be a Riemannian manifold and consider all $1$-arrows
on $M$ which preserve the metric $g.$ Such $1$-arrows define a $TLEFF$ $%
\mathcal{S}_{1}(M).$ We may fix some point $p\in M$ and fix some coordinates
around $p$ once and for all so that $g_{ij}(p)=\delta _{ij}$, thus
identifying $\mathcal{S}_{1}(M)_{p}^{p}$ with the orthogonal group $O(n).$
Now any $1$-arrow with source at $p$ defines an orthogonal frame at its
target $q$ by mapping the fixed orthogonal coordinate frame at $p$ to $q.$
The group $O(n)$ acts on all such $1$-arrows by composing on the source. Now
forgetting $1$-arrows but keeping the orthogonal frames defined by them, we
recover the orthogonal frame bundle of the metric $g$. However we will not
adapt this point of view. In view of the existence of geodesic coordinates,
we can now construct $2$-arrows on $M$ which preserve $1$-jet of $g$ , i.e., 
$1$-jet of $g$ at all $x\in M$ can be identified (in various ways). Thus we
obtain $\mathcal{S}_{2}(M)$ and the projection $\pi _{2,1}:\mathcal{S}%
_{2}(M)\rightarrow \mathcal{S}_{1}(M).$ As a remarkable fact, $\pi _{2,1}$
turns out to be an isomorphism. This fact is equivalent to the well known
Gauss trick of shifting the indices and showing the uniqueness of a metric
connection which is symmetric (Levi-Civita connection). The Christoffel
symbols are obtained now by twisting the $2$-arrows of $\mathcal{S}_{2}(M)$
by the $1$-arrows of $\mathcal{S}_{1}(M)$. Now we may not be able to
identify $2$-jet of $g$ over $M$ due to curvature of $g$ and thus we may
fail to construct the surjection $\pi _{3,2}:\mathcal{S}_{3}(M)\rightarrow 
\mathcal{S}_{2}(M).$ If we achive this (and $\pi _{3,2}$ will be again an
isomorphism), the next obstruction comes from the $3$-jet of $g$ which is
essentially the covariant derivative of the curvature. However, if $g$ has
constant curvature, then we can prolong $\mathcal{S}_{2}(M)$ uniquely to a
homogeneous geometry $\mathcal{S}_{\infty }(M)$ which, as a remarkable fact,
turns out to be a Klein geometry of order one since in this case $%
Lim_{k\rightarrow \infty }S_{k}(M)$ recaptures the isometry group of $(M,g)$
which acts transitively on $M$ and any isometry is uniquely determined by
any of its $1$-arrows$.$ Thus we may view a truncated geometry $\mathcal{S}%
_{k}(M)$ as a candidate for some homogeneous geometry $\mathcal{S}_{\infty
}(M)$ but $\mathcal{S}_{k}(M)$ must overcome the obstructions, if any, put
forward by $M.$ Almost all geometric structures (Riemannian, almost complex,
almost symplectic, almost quaternionic, ..) may be viewed as truncated
geometries of order at least one, each being a potential canditate for a
homogeneous geometry.

\begin{definition}
A truncated geometry $\mathcal{S}_{k}(M)$ is called formally integrable if
it prolongs to a homogeneous geometry.
\end{definition}

However, we require the prolongation required by Definition 6 to be
intrinsically determined by $\mathcal{S}_{k}(M)$ in some sense and not be
completely arbitrary. Given some $\mathcal{S}_{k}(M),$ note that Definition
6 requires the surjectivity of $\mathcal{S}_{j+1}(M)\rightarrow \mathcal{S}%
_{j}(M),$ $j\geqslant k.$ For instance, we may construct some $\mathcal{S}%
_{k+1}(M)$ in an intrinsic way without $\mathcal{S}_{k+1}(M)\rightarrow 
\mathcal{S}_{k}(M)$ being surjective. We may now redifine all lower terms by 
$\widetilde{\mathcal{S}_{j}(M)}=\pi _{k+1,j}\mathcal{S}_{k+1}(M)$ and start
anew at order $k+1.$ This is not allowed by Definition 6. For instance, let $%
\mathcal{S}_{1}(M)$ be defined by some almost symplectic form $\omega $ (or
almost complex structure $J).$ Then $\pi _{2,1}:\mathcal{S}%
_{2}(M)\rightarrow \mathcal{S}_{1}(M)$ will be surjective if $d\omega =0$ $($%
or $N(J)=0$ where $N(J)$ is the Nijenhuis tensor of $J$).

This prolongation process which we tried to sketch above is centered around
the concept of formal integrability which can be defined in full generality
turning the ambigious Definition 6 into a precise one. However, this
fundamental concept turns out to be highly technical and is fully developed
by D.C. Spencer and his co-workers from the point of view of PDEs,
culminating in [11]. More geometric aspects of this concept are emphasized
in [25], [26] and other books by this author .

\section{Bundle maps}

In this section we will briefly indicate the allowable bundle maps in the
present framework. Consider the universal $TLEFF$ $\mathcal{G}_{k}(M)$ of
order $k.$ We define the group bundle $\mathcal{AG}_{k}(M)\doteq \cup _{x\in
M}\mathcal{G}_{k}(M)_{x}^{x}$. The sections of this bundle form a group with
the operation defined fiberwise. We will call such a section a universal
bundle map (or a universal gauge transformation) of order $k.$ We will
denote the group of universal bundle maps by $\Gamma \mathcal{AG}_{k}(M).$
We obtain the projection $\pi _{k+1,k}:\Gamma \mathcal{AG}%
_{k+1}(M)\rightarrow \Gamma \mathcal{AG}_{k}(M)$ which is a homomorphism. If 
$\mathcal{S}_{k}(M)\subset \mathcal{G}_{k}(M),$ we similarly define $%
\mathcal{AS}_{k}(M)\subset \mathcal{AG}_{k}(M)$ and call elements of $\Gamma 
\mathcal{AS}_{k}(M)$ automorphisms (or gauge transformations) of $\mathcal{S}%
_{k}(M)$. Now let $\mathcal{S}_{k}(M)\subset \mathcal{G}_{k}(M$ and $%
g_{k}\in \Gamma \mathcal{AG}_{k}(M).$ We will denote $g_{k}(x)$ by $%
(g_{k})_{x}^{x},$ $x\in M.$ We define the $TLEFF$ $(Adg)\mathcal{S}_{k}(M)$
by defining its $k$-arrows as $(Adg)\mathcal{S}_{k}(M)_{y}^{x}\doteq
\{(g_{k})_{y}^{y}(f_{k})_{y}^{x}(g_{k}^{-1})_{x}^{x}\mid (f_{k})_{y}^{x}\in 
\mathcal{S}_{k}(M)_{y}^{x}\}.$

\begin{definition}
$\mathcal{S}_{k}(M)$ is called equivalent to $\mathcal{H}_{k}(M)$ if there
exists some $g\in \Gamma \mathcal{AG}_{k}(M)$ satisfying $(Adg)\mathcal{S}%
_{k}(M)=\mathcal{H}_{k}(M)$
\end{definition}

A necessary condition for the equivalence of $\mathcal{S}_{k}(M)$ and $%
\mathcal{H}_{k}(M)$ is that $\mathcal{S}_{k}(M)_{x}^{x}$ and $\mathcal{H}%
_{k}(M)_{x}^{x}$ be conjugate in $\mathcal{G}_{k}(M)_{x}^{x}$, i.e., they
must be compatible structures (like both Riemannian,...). Let $\frak{s}%
_{k}(M)$ and $\frak{h}_{k}(M)$ be the corresponding $TLEIF^{\prime }$s. The
above action induces an action of $g$ on $\frak{s}_{k}(M)$ which we will
denote also by $(Adg)\frak{s}_{k}(M).$ This latter action uses $1$-jet of $g$
since the geometric order of $\frak{s}_{k}(M)$ is $k+1,$ i.e., the
transformation rule of the elements (sections) of $\frak{s}_{k}(M)$ uses
derivatives up to order $k+1.$ This construction is functorial. In
particular, $(Adg)\mathcal{S}_{k}(M)=\mathcal{H}_{k}(M)$ implies $(Adg)\frak{%
s}_{k}(M)=\frak{h}_{k}(M)$. Clearly these actions commute with projections.
In this way we define the moduli spaces of geometries.

Some $g_{k}\in \Gamma \mathcal{AG}_{k}(M)$ acts on any $k^{th}$ order object
defined in Section 2. However, the action of $g_{k}$ does not commute with $%
d $ and therefore $g_{k}$ does not act on the horizontal complexes in (10).
To do this, we define $\mathcal{AG}_{\infty }(M)\doteq \cup _{x\in M}%
\mathcal{G}_{\infty }(M)_{x}^{x}$ where an element $(g_{\infty })_{x}^{x}$
of $\mathcal{G}_{\infty }(M)_{x}^{x}$ is the $\infty $-jet of some local
diffeomorphism with source and target at $x.$ As we noted above, $(g_{\infty
})_{x}^{x}$ is far from being a formal object: it determines this
diffeomorphism modulo the $\infty $-jet of identity diffeomorphism. Now some 
$g_{\infty }\in \Gamma \mathcal{AG}_{\infty }(M)$ does act on the horizontal
complexes in (10).

\section{Principal bundles}

Let $\mathcal{S}_{k}(M)\subset \mathcal{G}_{k}(M).$ We fix some $p\in M$ and
define $\mathcal{S}_{k}(M)^{(p)}\doteq \cup _{x\in M}\mathcal{S}%
_{k}(M)_{x}^{p}.$ The group $\mathcal{S}_{k}(M)_{p}^{p}$ acts on $\mathcal{S}%
_{k}(M)_{x}^{p}$ by composing with $k$-arrows of $\mathcal{S}_{k}(M)_{x}^{p}$
at the source as $(f_{k})_{x}^{p}\rightarrow (f_{k})_{x}^{p}(h_{k})_{p}^{p}$
and the projection $\mathcal{S}_{k}(M)^{(p)}\rightarrow M$ with fiber $%
\mathcal{S}_{k}(M)_{x}^{p}$ over $x$ is a principal bundle with group $%
\mathcal{S}_{k}(M)_{p}^{p}.$ Considering the adjoint action of $\mathcal{S}%
_{k}(M)_{p}^{p}$ on itself, we construct the associated bundle whose
sections are automorphisms (or gauge transformations) which we use in gauge
theory. This associated bundle can be identified with $\mathcal{AS}_{k}(M)$
in Section 6 and therefore the two concepts of automorphisms coincide. In
gauge theory, we let $h_{k}\in \Gamma \mathcal{AS}_{k}(M)$ act on $\mathcal{S%
}_{k}(M)^{(p)}$ on the target as $(f_{k})_{x}^{p}\rightarrow
(h_{k})_{x}^{x}(f_{k})_{x}^{p}$ reserving the source for the group $\mathcal{%
S}_{k}(M)_{p}^{p}.$ We will denote this transformation by $f\rightarrow
h\odot f.$ We can regard the object $h_{k}\odot \mathcal{S}_{1}(g)^{(p)}$ as
another principal $\mathcal{S}_{k}(M)_{p}^{p}$-bundle: we imagine two copies
of $\mathcal{S}_{k}(g)_{p}^{p}$, one belonging to principal bundle and one
outside which is the group of the principal bundle and $h_{k}$ acts only on
the principal bundle without changing the group. To be consistent with $%
\odot ,$ we now regard $h_{k}\in \Gamma \mathcal{AG}_{k}(M)$ as a general
bundle map and define the transform of $\mathcal{S}_{k}(M)^{(p)}$ by $h_{k}$
using $\odot .$ Now $\odot $ has a drawback from geometric point of view. To
see this, let $(M,g)$ be a Riemannian manifold. We will denote the $TLEFF$
determined by $g$ by $\mathcal{S}_{1}(g),$ identifying the principal $%
\mathcal{S}_{1}(g)^{(p)}\rightarrow M$ with the orthogonal frame bundle of $%
g $ and the group $\mathcal{S}_{1}(g)_{p}^{p}$ with $O(n)$ as in Section 5.
Now the transformed object $h\odot \mathcal{S}_{1}(g)^{(p)}$, which is
another $O(n)$-principal bundle, is not related to any metric in sight
unless $h$ = identity !! Thus we see that $\odot $ dispenses with the
concept of a metric but keeps the concept of an $O(n)$-principal bundle,
carrying us from our geometric envelope outside into the topological world
of general principal bundles. On the other hand, the action of $(\mathcal{G}%
_{1})_{x}^{x}$ on metrics at $x$ gives an action of $h$ on metrics on $M$
which we will denote by $g\rightarrow h\boxdot g.$ Changing our notation $%
(Adh)\mathcal{S}_{1}(g)$ defined in Section 6 to $h\boxdot \mathcal{S}%
_{1}(g) $ (using the same notation $\boxdot $), we see that $h\boxdot 
\mathcal{S}_{1}(g)=\mathcal{S}_{1}(h\boxdot g).$ Thus $\boxdot $ preserves
both metrics and also $1$-arrows determined by them.

Consider the naive inclusion

\begin{equation}
\text{differential geometry}\subset \text{topology}
\end{equation}

If we drop the word differential in (20), we may adapt the point of view
that the opposite inclusion holds now in (20). This is the point of view of
A. Grothendieck who views geometry as the study of \textit{form}, which
contains topology as a special branch (see his Promenade \#12, translated by
Roy Lisker). This broad perspective is clearly a far reaching generalization
of \textit{Erlangen Programm }presented here.

In view of (20), we believe that any theorem in the framework whose main
ingredients we attempted to outline here, however deep and far reaching, can
be formulated and proved also using principal bundles. To summarize, we may
say that differential geometry is the study \textit{smooth forms }and the
concepts of a \textit{form }and continuous deformation of \textit{forms }%
come afterwards as higher level of abstractions. We feel that it may be
fruitful to start with differential geometry rather than starting at a
higher level and then specializing to it.

\section{Connection and curvature}

Recall that a right principal $G$-bundle $P\rightarrow M$ determines the
groupoid $\frac{P\times P}{G}\rightarrow M\times M$ where the action of $G$
on $P\times P$ is given by $(x,y)g\doteq (xg,yg).$ Let $\mathcal{A}%
(P)\rightarrow M$ be the automorphism bundle obtained as the associated
bundle of $P\rightarrow M$ using the adjoint action of $G$ on itself as in
Section 7, whose sections are gauge transformations. We obtain in this way
the groupoid extension

\begin{equation}
1\longrightarrow \mathcal{A}(P)\longrightarrow \frac{P\times P}{G}%
\longrightarrow M\times M\longrightarrow 1
\end{equation}

where we again use the first and last arrows to indicate injectivity and
surjectivity without any algebraic meaning. On the infinitesimal level, (21)
gives the Atiyah sequence of $P\rightarrow M$

\begin{equation}
0\longrightarrow \mathcal{LA}(P)\longrightarrow \frac{TP}{G}\overset{\pi }{%
\longrightarrow }TM\longrightarrow 0
\end{equation}

(see [19], [20] for the details of the Atiyah sequence) where $\mathcal{LA}%
(P)\longrightarrow M$ is the Lie algebra bundle obtained as the associated
bundle using the adjoint action of $G$ on its Lie algebra $\mathcal{L}(G).$
Connection forms $\omega $ on $P\rightarrow M$ are in 1-1 correspondence
with transversals in (22), i.e., vector bundle maps $\omega :TM\rightarrow 
\frac{TP}{G}$ with $\pi \circ \omega =id$ and curvature $\kappa $ of $\omega 
$ is defined by $\kappa (X,Y)=\omega \lbrack X,Y]-[\omega X,\omega Y].$ The
extension (22) splits iff (22) admits a flat transversal. Thus Atiyah
sequence completely recovers the formalism of connection and curvature on $%
P\rightarrow M$ in the framework of algebroid extensions as long as we work
over a fixed base manifold $M$.

Now let $\mathcal{S}_{\infty }(M)$ be a homogeneous geometry with
infinitesimal geometry $\frak{s}_{\infty }(M).$ The groupoid $\frac{\mathcal{%
S}_{k}(M)^{(p)}\times \mathcal{S}_{k}(M)^{(p)}}{\mathcal{S}_{k}(M)_{p}^{p}}$
in (21) determined by the principal bundle $\mathcal{S}_{k}(M)^{(p)}%
\rightarrow M$ as defined in Section 7 can be identified with $\mathcal{S}%
_{k}(M)$ and the algebroid $\frac{T\mathcal{S}_{k}(M)^{(p)}}{\mathcal{S}%
_{k}(M)_{p}^{p}}$ in (22) can be identified with $\frak{s}_{k}(M).$ We have $%
\mathcal{A}(\mathcal{S}_{k}(M)^{(p)})=\cup _{x\in M}\mathcal{S}%
_{k}(M)_{x}^{x}$ as already indicated in Section 7 and $\mathcal{LA}(%
\mathcal{S}_{k}(M)^{(p)})$ $\doteq \cup _{x\in M}\mathcal{L}(\mathcal{S}%
_{k}(M)_{x}^{x})$ where the bracket of sections is defined fiberwise. Thus
(21) becomes

\begin{equation}
1\longrightarrow \mathcal{AS}_{k}(M)\longrightarrow \mathcal{S}%
_{k}(M)\longrightarrow M\times M\longrightarrow 1
\end{equation}

and the Atiyah sequence (22) is now

\begin{equation}
0\longrightarrow \mathcal{LAS}_{k}(M)\longrightarrow \frak{s}%
_{k}(M)\longrightarrow TM\longrightarrow 0
\end{equation}

It is easy to check exactness of (24) in local coordinates using (2).

Our purpose is now to indicate how the present framework captures
information peculiar to jets by changing the base manifold, which Atiyah
sequence does not detect. To see this, let $m\leq k+1$ and consider the Lie
group extension at $x\in M:$

\begin{equation}
1\longrightarrow \mathcal{S}_{k,m}(M)_{x}^{x}\longrightarrow \mathcal{S}%
_{k}(M)_{x}^{x}\longrightarrow \mathcal{S}_{m}(M)_{x}^{x}\longrightarrow 1
\end{equation}

where the kernel $\mathcal{S}_{k,m}(M)_{x}^{x}$ is nilpotent if $m\geq 1$
and is abelian if $k=m+1\geq 2.$ Consider the $\mathcal{S}_{k,m}(M)_{p}^{p}$%
-principal bundle $\mathcal{S}_{k}(M)^{(p)}\rightarrow \mathcal{S}%
_{m}(M)^{(p)}.$ If $\mathcal{S}_{k,m}(M)_{p}^{p}$ is contractible (this is
the case in many examples for $m\geq 1)$, this principal bundle is trivial
and its Atiyah sequence admits flat transversals. Thus nothing is gained by
considering higher order jets and all information is contained in the Atiyah
sequence of $\mathcal{S}_{1}(M)^{(p)}\rightarrow M.$ On the other hand, we
have the following extension of $TLEIF^{\prime }s:$%
\begin{equation}
0\longrightarrow \mathcal{LS}_{k,m}(M)\longrightarrow \frak{s}%
_{k}(M)\longrightarrow \frak{s}_{m}(M)\longrightarrow 0
\end{equation}

where $\mathcal{LS}_{k,m}(M)\doteq \cup _{x\in M}\mathcal{L}(\mathcal{S}%
_{k,m}(M)_{x}^{x}).$ Using (2), it is easy to check that the existence of a
flat transversal in (26) a priori forces the splitting of the Lie algebra
extension

\begin{equation}
0\longrightarrow \mathcal{L}(\mathcal{S}_{k,m}(M)_{x}^{x})\longrightarrow 
\mathcal{L}(\mathcal{S}_{k}(M)_{x}^{x})\longrightarrow \mathcal{L}(\mathcal{S%
}_{m}(M)_{x}^{x})\longrightarrow 0
\end{equation}

for all $x\in M$ where (27) is (25) in infinitesimal form. However (25) and
(27) do not split in general. For instance, (27) does not split even in the
universal situation $\mathcal{S}_{\infty }(M)=\mathcal{G}_{\infty }(M)$ when 
$m=2$ and $k=3.$ In fact, the dimensions of the extension groups $H^{2}(%
\mathcal{L}(\mathcal{S}_{m}(M)_{x}^{x}),\mathcal{S}_{m+1,m}(M)_{x}^{x})$ are
computed in [28] for all $m$ when $\dim M=1.$

Thus we see that the theory of principal bundles, which is essentially
topological, concentrates on the maximal compact subgroup of the structure
group of the principal bundle as it is this group which produces nontrivial
characteristic classes as invariants of equivalence classes of principal
bundles modulo bundle maps. Consequently, this theory concentrates on the
contractibility of the kernel in (25) whereas it is the types of the
extensions in (25), (27) which emerge as the new ingredient in the present
framework.

The connections considered above are transversals and involve only $%
TLEIF^{\prime }s$ (algebroids). There is another notion of connection based
on Maurer-Cartan form, which is incorporated by the nonlinear Spencer
sequence (see Theorem 31 on page 224 in [25]), where the passage from $TLEFF$
(groupoid) to its $TLEIF$ (algebroid) is used in a crucial way. The passage
from extensions of $TLEFF^{\prime }s$ (torsionfree connections in finite
form) to extensions of their $TLEIF^{\prime }s$ (torsionfree connections in
infinitesimal form) relates these two notions by means of a single diagram
which we hope to discuss elsewhere. The approach to parabolic geometries
adapted in [5] is a complicated mixture of these two notions whose
intricacies, we believe, will be clearly depicted by this diagram.

However we view a connection, the main point here seems to be that it
belongs to the group rather than to the space on which the group acts. Since
there is an abundance of groups acting transitively on some given space, it
seems meaningles to speak of the curvature of some space unless we specify
the group. However, it turns out that the knowledge of the $k$-arrows of
some ideal group is sufficient to define a connection but this connection
will not be unique except in some special cases.

\section{Some remarks}

In this section (unfortunately somewhat long) we would like to make some
remarks on the relations between the present framework and some other
frameworks.

$1)$ Let $\mathcal{E}\rightarrow M$ be a differentiable fibered space. It
turns out that we have an exterior calculus on $J^{\infty }(\mathcal{E}%
)\rightarrow M.$ Decomposing exterior derivative and forms into their
horizontal and vertical components, we obtain a spectral sequence, called
Vinogradov $\mathcal{C}$-spectral sequence, which is fundamental in the
study of calculus of variations (see [32], [33] and the references therein).
The limit term of $\mathcal{C}$-spectral sequence is $H_{deR}^{\ast
}(J^{\infty }(\mathcal{E})).$ In particular, if $\mathcal{E}=T(M),$ we
obtain $H_{deR}^{\ast }(\frak{g}_{\infty }(M)).$ The $\mathcal{C}$-spectral
sequence can be defined also with coefficients ([21]). On the other hand, we
defined $H^{\ast }(\frak{g}_{\infty }(M),J_{\infty }(M))$ and $H^{\ast }(%
\frak{g}_{\infty }(M),\frak{s}_{\infty }(M),J_{\infty }(M))$ in Sections 2,
3. As we indicated in Section 2, we can consider representations of $\frak{g}%
_{\infty }(M)$ other than $J_{\infty }(M)$ (for instance, see $2)$ below).
These facts hint, we feel, at the existence of a very general Van Est type
theorem which relates these cohomology groups.

$2)$ Recall that (5) and (7) define a representation of the algebroid $\frak{%
g}_{k}(M)$ on the vector bundle $J_{k}(M)\rightarrow M$. Thus we can
consider the cohomology groups $H^{\ast }(\frak{g}_{k}(M),J_{k}(M)$ as
defined in [19], [20] (see also the references there for original sources)
and $H^{\ast }(\frak{g}_{k}(M),J_{k}(M))$ coincides with the cohomology of
the bottom sequence of (17) by definition. Now $\frak{g}_{k}(M)$ has other
``intrinsic'' representations, for instance $\frak{g}_{k-1}(M).$ Lemma 8.32
and the formula on page 383 in [25] (which looks very similar to (4)) define
this representation. In particular, the cohomology groups $H^{\ast }(\frak{g}%
_{k}(M),$ $\frak{g}_{k-1}(M))$ are defined and are given by the bottom
sequence of (17) for this case. Using this representation of $\frak{g}%
_{k}(M) $ on $\frak{g}_{k-1}(M)$, deformations of $TLEIF^{\prime }$s are
studied in [25] in detail using Janet sequence (Chapter 7, Section 8 of
[25]). Recently, some deformation cohomology groups are introduced in [8] in
the general framework of algebroids. However, if the algebroid is a $TLEIF,$
it seems to us that these cohomology groups coincide with those in [25] (and
therefore also with the bottom sequence of (17)) and are not new (see also
[20], pg. 309 for a similar claim). In view of the last paragraph of Section
4, the fact that deformation cohomology arises as sheaf cohomology in
Kodaira-Spencer theory and as algebroid cohomology in [25], [8] is no
coincidence.

$3)$ Let $\mathcal{S}_{2}(M)$ be a truncated geometry on $M.$ We fix some $%
x\in M$ and consider the following diagram of Lie group extensions

\begin{equation}
\begin{array}{ccccccccc}
0 & \longrightarrow & \mathcal{G}_{2,1}(M)_{x}^{x} & \longrightarrow & 
\mathcal{G}_{2}(M)_{x}^{x} & \longrightarrow & \mathcal{G}_{1}(M)_{x}^{x} & 
\longrightarrow & 1 \\ 
&  & \uparrow &  & \uparrow &  & \uparrow &  &  \\ 
0 & \longrightarrow & \mathcal{S}_{2,1}(M)_{x}^{x} & \longrightarrow & 
\mathcal{S}_{2}(M)_{x}^{x} & \longrightarrow & \mathcal{S}_{1}(M)_{x}^{x} & 
\longrightarrow & 1
\end{array}
\end{equation}

where the vertical maps are inclusions. The top row of (28) splits and the
components of these splittings are naturally interpreted as the Christoffel
symbols of symmetric ``point connections''. Some $k_{x}^{x}\in \mathcal{G}%
_{2,1}(M)_{x}^{x}$, which is a particular bundle map as defined in Section 6
when $M=\{x\},$ transforms $\mathcal{S}_{2}(M)_{x}^{x}$ by conjugation but
acts as identity on $\mathcal{S}_{2,1}(M)_{x}^{x}$ and $\mathcal{S}%
_{1}(M)_{x}^{x}.$ Using $\mathcal{G}_{2,1}(M)_{x}^{x}$ as allowable
isomorphisms, we defined in [13] the group of \textit{restricted }extensions 
$H_{res}^{2}(\mathcal{S}_{1}(M)_{x}^{x},\mathcal{S}_{2,1}(M)_{x}^{x})$. This
group vanishes iff the restriction of some splitting of the top row of (28)
to the bottom row splits also the bottom row or equivalently, iff the bottom
row admits a symmetric point connection. The main point is that this group
is sensitive to phenomena happening only inside the top row of (28) which is
our universal envelope as in Section 2. Using the Lie algebra anolog of
(28), we defined also the group $H_{res}^{2}(\mathcal{L}(\mathcal{S}%
_{1}(M)_{x}^{x}),\mathcal{S}_{2,1}(M)_{x}^{x})$ (note that $\mathcal{S}%
_{2,1}(M)_{x}^{x}\subset \mathcal{G}_{2,1}(M)_{x}^{x}$ are vector groups)
obtaining the homomorphism $H_{res}^{2}(\mathcal{S}_{1}(M)_{x}^{x},\mathcal{S%
}_{2,1}(M)_{x}^{x})\rightarrow H_{res}^{2}(\mathcal{L}(\mathcal{S}%
_{1}(M)_{x}^{x}),\mathcal{S}_{2,1}(M)_{x}^{x})$ ([13]). On the other hand,
regarding the bottom row of (28) as an \textit{arbitrary }Lie group
extension as in [15] without any reference to our universal envelope, we can
define $H^{2}(\mathcal{S}_{1}(M)_{x}^{x},\mathcal{S}_{2,1}(M)_{x}^{x})\ $\
and the homomorphism $H^{2}(\mathcal{S}_{1}(M)_{x}^{x},\mathcal{S}%
_{2,1}(M)_{x}^{x})\rightarrow H^{2}(\mathcal{L}(\mathcal{S}_{1}(M)_{x}^{x}),%
\mathcal{S}_{2,1}(M)_{x}^{x})$. Thus we obtain the following commutative
diagram

\begin{equation}
\begin{array}{ccc}
H_{res}^{2}(\mathcal{S}_{1}(M)_{x}^{x},\mathcal{S}_{2,1}(M)_{x}^{x}) & 
\longrightarrow & H_{res}^{2}(\mathcal{L}(\mathcal{S}_{1}(M)_{x}^{x}),%
\mathcal{S}_{2,1}(M)_{x}^{x}) \\ 
\downarrow &  & \downarrow \\ 
H^{2}(\mathcal{S}_{1}(M)_{x}^{x},\mathcal{S}_{2,1}(M)_{x}^{x} & 
\longrightarrow & H^{2}(\mathcal{L}(\mathcal{S}_{1}(M)_{x}^{x}),\mathcal{S}%
_{2,1}(M)_{x}^{x})
\end{array}
\end{equation}

where the vertical homomorphisms are induced by inclusion. In [13], we gave
examples where $H_{res}^{2}(\mathcal{L}(\mathcal{S}_{1}(M)_{x}^{x}),\mathcal{%
S}_{2,1}(M)_{x}^{x})$ is nontrivial whereas $H^{2}(\mathcal{L}(\mathcal{S}%
_{1}(M)_{x}^{x}),$

$\mathcal{S}_{2,1}(M)_{x}^{x})$ is trivial. This fact shows that we may
loose information when we pass from the top row to the bottom row in (29).

Now the action of $\mathcal{S}_{1}(M)_{x}^{x}$ on $\mathcal{S}%
_{2,1}(M)_{x}^{x}$ gives a representation of the $TLEFF$ $\mathcal{S}_{1}(M)$
on the vector bundle $\mathcal{S}_{2,1}(M)\doteq \cup _{x\in M}\mathcal{S}%
_{2,1}(M)_{x}^{x}$ and thus we can define the cohomology groups $%
H_{res}^{\ast }(\mathcal{S}_{1}(M),\mathcal{S}_{2,1}(M))$ in such a way that
they will respect our universal envelope. In this way we arrive at the
global analog of (29). We believe that we will loose information also in
this global case. To summarize, even though the constructions in this paper
can be formulated in the general framework of groupoids and algebroids as in
[7], [8],[19], [20], we believe that this general framework will not be
sensitive in general to certain phenomena peculiar to jets unless it takes
the universal homogeneous envelope into account. In particular, we would
like to express here our belief that Lie equations form the geometric core
of groupoids which are the ultimate generalizations of (pseudo)group actions
in which we dispense with the action but retain the symmetry that the action
induces on the space (compare to the introduction of [36]).

$4)$ Let $\mathcal{X}$ $(M)$ be the Lie algebra of smooth vector fields on $%
M.$ Recalling that $j_{k}(\mathcal{X}(M))\subset \frak{g}_{k}(M)$, (5) gives
a representation of $\mathcal{X}(M)$ on $J_{k}(M)$. Denoting the cochains
computing this cohomology by $\overset{(k,r)}{\wedge }_{GF},$ we obtain the
chain map $\overset{(k,\ast )}{\wedge }\rightarrow $ $\overset{(k,\ast )}{%
\wedge }_{GF}$ which indicates that Gelfand-Fuks cohomology is involved in
the present framework and plays a central role.

$5)$ This remark can be considered as the continuation of Section 7 and the
last paragraph of Section 4.

Let $Q$ be the subgroup of $\mathcal{G}_{1}(n+1)=GL(n+1,\mathbb{R)}$
consisting of matrices of the form

\begin{equation}
\left[ 
\begin{array}{cc}
A & 0 \\ 
\xi & \lambda
\end{array}
\right]
\end{equation}
where $A$ is an invertible $n\times n$ matrix, $\xi =(\xi _{1},...,\xi _{n})$
and $\lambda \neq 0.$ We will denote (30) by $(A,\xi ,\lambda ).$ We have
the homomorphism $Q/\lambda I\rightarrow \mathcal{G}_{1}(n)$ defined by $%
(A,\xi ,1)\rightarrow A,$ where $\lambda I$ denotes the subgroup $\{\lambda
I\mid \lambda \in \mathbb{R\}},$ with the abelian kernel $K$ consisting of
elements of the form $(I,\xi ,1).$ We also have the injective homomorphism $%
Q/\lambda I\rightarrow \mathcal{G}_{2}(n)$ defined by $(A_{j}^{i},\xi
_{j},1)\rightarrow (A_{j}^{i},\xi _{j}A_{k}^{i}+\xi _{k}A_{j}^{i})$ which
gives the diagram

\begin{equation}
\begin{array}{ccccccccc}
0 & \longrightarrow & \mathcal{G}_{2,1}(n) & \longrightarrow & \mathcal{G}%
_{2}(n) & \longrightarrow & \mathcal{G}_{1}(n) & \longrightarrow & 1 \\ 
&  & \uparrow &  & \uparrow &  & \parallel &  &  \\ 
0 & \longrightarrow & K & \longrightarrow & Q/\lambda I & \longrightarrow & 
\mathcal{G}_{1}(n) & \longrightarrow & 1
\end{array}
\end{equation}

Now $(\mathcal{G}_{1}(n+1),Q)$ is a Klein pair of order two with ghost $N=$ $%
\lambda I$. We also have the effective Klein pair $(\mathcal{G}%
_{1}(n+1)/\lambda I,Q/\lambda I)$ which is of order two. Note that the
standard action of $\mathcal{G}_{1}(n+1)$ on $\mathbb{R}^{n+1}\backslash 0$
induces a transitive action of $\mathcal{G}_{1}(n+1)$ on the the real
projective space $\mathbb{R}P(n)$ and $Q$ is the stabilizer of the coloumn
vector $p=$ $(0,0,...,1)^{T}$ so that both Klein pairs $(\mathcal{G}%
_{1}(n+1),Q)$ and $(\mathcal{G}_{1}(n+1)/\lambda I,Q/\lambda I)$ define the
same base $\mathbb{R}P(n).$ Clearly, $\mathcal{G}_{1}(n+1)$ and $\mathcal{G}%
_{1}(n+1)/\lambda I$ induce the same $2$-arrows on $\mathbb{R}P(n).$

At this stage, we have two relevant principle bundles.

$i)$ The principle bundle $\mathcal{G}_{1}(n+1)^{(p)}\rightarrow \mathbb{R}%
P(n)$ which consists of all $2$-arrows eminating from $p$ and has the
structure group $Q/\lambda I.$ This is the same principle bundle as $(%
\mathcal{G}_{1}(n+1)/\lambda I)^{(p)}\rightarrow \mathbb{R}P(n)$. In
particular, (18) reduces to (31) in this case (the reader may refer to
Example 4.1 on pg. 132 of \ [18] and also to the diagram on pg.142).

$ii)$ The principle bundle $\mathcal{G}_{1}(n+1)\rightarrow \mathbb{R}P(n)$
with structure group $Q.$ Note that we have the central extension

\begin{equation}
1\longrightarrow \lambda I\longrightarrow Q\longrightarrow Q/\lambda
I\longrightarrow 1
\end{equation}

which splits: some $q=(A,\xi ,\lambda )\in Q$ factors as $q=ab$ where $%
a=(\lambda ^{-1}A,\lambda ^{-1}\xi ,1)$

$\in \mathcal{G}_{2}(n)$ and $b=\lambda I.$ Note that $ii)$ is obtained from 
$i)$ by lifting the structure group $Q/\lambda I$ to $Q$ in (31).

Now we have a representation of $Q$ on $\mathbb{R}$ defined by the
homomorphism $(A,\xi ,\lambda )\rightarrow \lambda ^{-N}$ for some integer $%
N\geq 0.$ Replacing $\mathbb{R}$ by $\mathbb{C}$ and working with complex
groups and holomorphic actions, it is known that the holomorphic sections of
the associated line bundle of $\mathcal{G}_{1}(n+1,\mathbb{C})\rightarrow 
\mathbb{C}P(n)$ realize all irreducable representations of the unitary group 
$U(n)$ as $N$ varies (see [16], pg. 138-152 and [17]). As a very crucial
fact, we can repeat this construction by replacing the Klein pair $(\mathcal{%
G}_{1}(n+1,\mathbb{C}),Q)$ of order two by the effective Klein pair $%
(U(n),U(n-1)\times U(1))$ of order one and this latter construction recovers
the same line bundle (see [16] for details). The following question
therefore arises naturally: Let $(G,H)$ be a Klein pair with ghost $N.$ Let $%
\rho :H\rightarrow GL(V)$ be a representation and $E\rightarrow M=G/H$ be
the associated homogeneous vector bundle of $G\rightarrow G/H.$ Can we
always find some \textit{effective} Klein pair $(\overline{G},\overline{H})$
(not necessarily of the same order) with $\overline{G}/\overline{H}=M,$ a
representation $\overline{\rho }:\overline{H}\rightarrow GL(V)$ such that $%
E\rightarrow M$ is associated with $(\overline{G})^{(p)}\rightarrow M$, or
shortly

$\mathbf{Q2:}$ Can we always avoid ghosts in Klein geometry?

Replacing $\lambda I,$ $Q,$ $Q/\lambda I$ in (31) respectively by $U(1),$ $%
Spin^{c}(4),$ $SO(4)$ and recalling the construction of $Spin^{c}$-bundle on
a 4-manifold ([22]), we see that $\mathbf{Q2}$ is quite relevant as it asks
essentially the scope and limitations of \textit{Erlangen Programm.}

$6)$ The assumption of transitivity, i.e., the surjectivity of the right
arrows in (13) , (14), is imposed upon us by \textit{Erlangen Programm.}
However, many of the constructions in this paper can be carried out without
the assumption of transitivity. For instance, foliations give rise to
intransitive Lie equations but they are studied in the literature mostly
from the point of view of general groupoids and algebroids (see the
references in [7], [8]).

Our main object of study in this paper has been a differentiable manifold $%
M. $ The sole reason for this is that this author has been obsessed years
ago by the question ``What are Christoffel symbols? '' and he could not
learn algebraic geometry from books and he did not have the chance to learn
it from experts (this last remark applies also to differential geometry) as
he has always been at the wrong place at the right time. We feel (and
sometimes almost see, for instance [27], [3]) that the present framework has
also an algebraic counterpart valid for algebraic varieties.

\bigskip

\bigskip

\bigskip

\bigskip

\bigskip

\textbf{References}

\bigskip

[1] E.Abado\u{g}lu, preprint

[2] A.Banyaga, The structure of Classical Diffeomorphism Groups, Kluwer
Academic Publishers, Volume 400, 1997

[3] A.Beilinson, V.Ginzburg: Infinitesimal structure of moduli spaces of
G-bundles, Internat. Math. Res. Notices, no 4, 93-106, 1993

[4] R.Bott: Homogeneous vector bundles, Ann. of Math., Vol. 66, No. 2,
203-248, 1957

[5] A.Cap, J.Slovak, V.Soucek: Bernstein-Gelfand-Gelfand sequences, Ann. of
Math. 154, 97-113, 2001

[6] C.Chevalley, S.Eilenberg: Cohomology theory of Lie groups and Lie
algebras. Trans. Amer. Math. Soc. 63, (1948). 85--124

[7] M.Crainic: Differentiable and algebroid cohomology, Van Est isomorphism,
and characteristic classes, Comment. Math. Helv. 78, 2003, 681-721

[8] M.Crainic, I.Moerdijk: Deformations of Lie brackets: cohomological
aspects, arXiv: 0403434

[9] W.T.Van Est: Group cohomology and Lie Algebra cohomology in Lie groups,
Indagationes Math., 15, 484-492, 1953

[10] D.B.Fuks: Cohomology of infinite-dimensional Lie algebras. Translated
from the Russian by A.B. Sosinski, Contemporary Soviet Mathematics,
Consultants Bureau, New York, 1986

[11] H.Goldschmidt, D.C.Spencer: On the nonlinear cohomology of Lie
equations, I, II, III, IV, Acta Math. 136, 103-170, 171-239, 1976, J.
Differential Geometry, 13, 409-453, 455-526, 1979

[12] W.Greub, S.Halperin, R.Vanstone: Connections, Curvature and Cohomology,
Vol.III, Academic Press, New York San Fransisco London, 1976

[13] B.G\"{u}rel, E.Orta\c{c}gil, F.\"{O}zt\"{u}rk: Group extensions in
second order jet group, preprint

[14] T.Hawkins: Emergence of the Theory of Lie Groups, An Essay in the
History of Mathematics 1869-1926, Sources and Studies in the History of
Mathematics and Physical Sciences, 2000 Springer-Verlag New York, Inc.

[15] G.Hochschild: Group extensions of Lie groups. Ann. of Math. (2) 54,
(1951). 96--109

[16] A.W.Knapp: Lie Groups, Lie Algebras and Cohomology, Mathematical Notes
34, Princeton University Press, NJ, 1988

[17] A.W.Knapp: Representation Theory of Semisimple Groups: An Overview
based on Examples, Mathematical Notes, Princeton University Press,
Princeton, NJ, 1986

[18] S.Kobayashi: Transformation Groups in Differential Geometry, Ergebnisse
der Mathematic und ihrer Grenzgebiete. Band 70, Springer-Verlag, 1972

[19] K.Mackenzie: Lie Groupoids and Lie Algebroids in Differential Geometry,
London Mathematical Society Lecture Note Series, 124, Cambridge University
Press, Cambridge, 1987

[20] K.Mackenzie: General Theory of Lie Groupoids and Lie Algebroids, London
Mathematical Society Lecture Note Series, 213, Cambridge University Press,
Cambridge, 2005

[21] M. Marvan: On the $\mathcal{C}$-spectral sequence with ''general''
coefficients. Differential geometry and its applications (Brno, 1989),
361--371, World Sci. Publishing, Teaneck, NJ, 1990

[22] J.W.Morgan: The Seiberg-Witten Equations and Applications to the
topology of smooth Four-Manifolds, Mathematical Notes 44, Princeton
Univertsity Press, Princeton, NJ, 1996

[23] H.Omori: Infinite dimensional Lie transformation groups. Lecture Notes
in Mathematics, Vol. 427. Springer-Verlag, Berlin-New York, 1974

[24] R.S.Palais: Extending diffeomorphisms, Proc. Amer. Math. Soc. 11 1960
274--277

[25] J.F. Pommaret: Systems of Partial Differential Equations and Lie
Pseudogroups, Gordon and Breach Science Publishers, New York, London, Paris,
1978

[26] J.F.Pommaret: Partial Differential Equations and Group Theory. New
perspectives for applications. Mathematics and its Applications, 293. Kluwer
Academic Publishers Group, Dordrecht, 1994

[27] Z.Ran: Derivatives of Moduli, Internat. Math. Res. Notices, no 4,
63-74, 1992

[28] B.K.Reinhart: Some remarks on the structure of the Lie algebra of
formal vector fields. Transversal structure of foliations (Toulouse, 1982).
Ast\'{e}risque No. 116 (1984), 190--194

[29] R.W.Sharpe: Differential Geometry, Cartan's Generalization of Klein's
Erlangen Program, Graduate Texts in Mathematics, Springer-Verlag, New York
Berlin Heidelberg, 1997

[30] W.Thurston: Three-dimensional Geometry and Topology Vol 1, Edited by
Silvio Levy, Princeton Mathematical Series, 35, Princeton University Press,
NJ, 1997

[31] G.Vezzosi, A.M.Vinogradov: On higher order analogues of de Rham
cohomology, Differential Geom. Appl. 19 (2003), no. 1, 29--59

[32] A.M.Vinogradov: Cohomological Analysis of Partial Differential
Equations and Secondary Calculus, Translations of Mathematical Monographs,
Volume 204, AMS, 2000

[33] A.M.Vinogradov: Scalar differential invariants, diffieties and
characteristic classes, Mechanics, Analysis and Geometry, 200 years after
Lagrange, M.Francaviglia (editor), Elsevier Science Publishers B.V., 1991,
379-414

[34] H.C.Wang: Closed manifolds with homogeneous structures, Amer.J.Math.,
76, 1954, 1-32

[35] G.Weingart: Holonomic and semi-holonomic geometries, Seminaires \&
Congres, 4, 2000, 307-328

\bigskip \lbrack 36] A.Weinstein: Groupoids: unifying internal and external
symmetry. A tour through some examples. Groupoids in analysis, geometry, and
physics (Boulder, CO, 1999), 1--19, Contemp. Math., 282, Amer. Math. Soc.,
Providence, RI, 2001

\bigskip

\bigskip

Erc\"{u}ment Orta\c{c}gil, Bo\u{g}azi\c{c}i University, Bebek, 34342,
Istanbul, Turkey

e-mail: ortacgil@boun.edu.tr

\end{document}